\documentclass[12pt, a4paper]{article}
\usepackage{graphicx, theorem, color} \definecolor{webblue}{rgb}{0,0,0.6}
\usepackage[colorlinks=true, linkcolor=webblue, citecolor=webblue,
 urlcolor=webblue, pdfstartview=FitBH,
 pdftitle={Quadratic matings and ray connections},
 pdfauthor={Wolf Jung},
 bookmarksopen=true, bookmarksnumbered=true, pdfpagemode=UseNone]{hyperref}
\makeatletter 
\long\def\@makecaption#1#2{%
  \vskip\abovecaptionskip
  {\small\textbf{#1:} #2}\par
  \vskip\belowcaptionskip}
\makeatother
\theoremstyle{break}
\theorembodyfont{\itshape}  \newtheorem{thm}{Theorem}[section]
\newtheorem{dfn}[thm]{Definition} \newtheorem{prop}[thm]{Proposition}
 \newtheorem{conj}[thm]{Conjecture}

\theorembodyfont{\slshape}  \newtheorem{rmk}[thm]{Remark}  
\newtheorem{xmp}[thm]{Example}   
\newcommand{\be}{\begin{equation}} \newcommand{\ee}{\end{equation}}
\newcommand{\ban}{\begin{eqnarray}} \newcommand{\ean}{\end{eqnarray}}
\renewcommand{\hat}{\widehat} \renewcommand{\tilde}{\widetilde}
\newcommand{\mybox}{\hspace*{\fill}\rule{2mm}{2mm}}
\DeclareSymbolFont{AMSb}{U}{msb}{m}{n}
\DeclareSymbolFontAlphabet{\mathbb}{AMSb}
\newcommand{\C}{\mathbb{C}}
  \newcommand{\R}{\mathbb{R}}
\newcommand{\Z}{\mathbb{Z}}
  \newcommand{\disk}{\mathbb{D}}

\let\inodot\i  \renewcommand{\i}{\mathrm{i}}
\newcommand\ddoti{\"\inodot}  \newcommand\brevei{\u\inodot}
\renewcommand{\phi}{\varphi}
\renewcommand{\r}{\mathcal{R}}
\newcommand{\J}{\mathcal{J}}  \newcommand{\K}{\mathcal{K}}
\newcommand{\M}{\mathcal{M}}  \newcommand{\sM}{{\scriptscriptstyle M}}

\newcommand{\eqa}{\simeq} 
\newcommand{\eqg}{\cong} 
\newcommand{\fmate}{\sqcup} 
\newcommand{\tmate}{\coprod} 
\newcommand{\2}{$(2,\,2,\,2,\,2)$}
\newcommand{\4}{$(2,\,4,\,4)$}
\date{}
\author{Wolf Jung\\
{\small Gesamtschule Brand, 52078 Aachen, Germany,}\\
{\small and Jacobs University, 28759 Bremen, Germany.}\\
{\small E-mail: \href{mailto:jung@mndynamics.com}{jung@mndynamics.com}}}
\title{Quadratic matings and ray connections}
\begin{document}
\maketitle
\begin{abstract}\noindent
A topological mating is a map defined by gluing together the filled Julia
sets of two quadratic polynomials.  The identifications are visualized and
understood by pinching ray-equivalence classes of the formal mating.  For
postcritically finite polynomials in non-conjugate limbs of the Mandelbrot
set,  classical results construct the geometric mating from the formal
mating.  Here families of examples are discussed,  such that all
ray-equivalence classes are uniformly bounded trees.  Thus the topological
mating is obtained directly in geometrically finite and infinite cases.  On
the other hand,  renormalization provides examples of unbounded cyclic ray
connections,  such that the topological mating is not defined on a
Hausdorff space.
\par\noindent
There is an alternative construction of mating,  when at least one polynomial
is preperiodic:  shift the infinite critical value of the other polynomial to
a preperiodic point.  Taking homotopic rays, it gives simple examples of
shared matings.  Sequences with unbounded multiplicity of sharing,  and
slowly growing preperiod and period,  are obtained both in the Chebychev
family and for Airplane matings.  Using preperiodic polynomials with
identifications between the two critical orbits,  an example of mating
discontinuity is described as well.
\end{abstract}

\section{Introduction} \label{1}
Starting from two quadratic polynomials $P(z)=z^2+p$ and $Q(z)=z^2+q$,
construct the topological mating $P\tmate Q$ by gluing the filled Julia sets
$\K_p$ and $\K_q$\,. If there is a conjugate rational map $f$, this defines
the geometric mating. These maps are understood by starting with the formal
mating $g=P\fmate Q$, which is conjugate to $P$ on the lower half-sphere
$|z|<1$ and to $Q$ on the upper half-sphere $|z|>1$ of
$\hat\C=\C\cup\{\infty\}$\,: ray-equivalence classes consist of external rays
of $P$ and $Q$ with complex conjugate angles, together with landing points in
$\partial\K_p$ and $\partial\K_q$\,; collapsing these classes defines the
topological mating. In the postcritically finite case, with $p$ and $q$ not in
conjugate limbs of $\M$, either $g$ or a modified version
$\tilde g$ is combinatorially equivalent and semi-conjugate to a rational map
$f$ \cite{rst, teich, dh, book2h, rs}. So the topological mating exists and
$f$ is conjugate to it --- it is a geometric mating.

\begin{figure}[h!t!b!]
\unitlength 0.001\textwidth 
\begin{picture}(990, 420)
\put(10, 0){\includegraphics[width=0.42\textwidth]{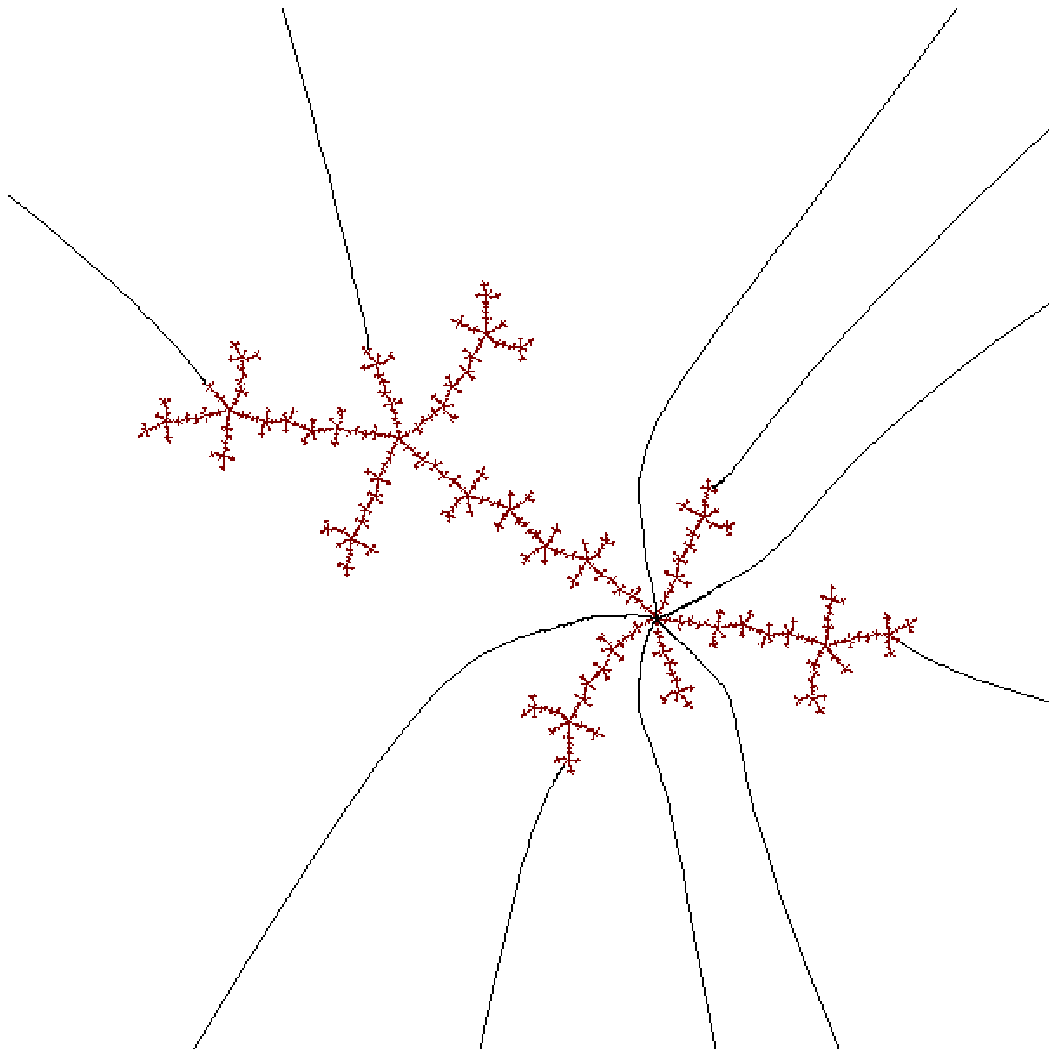}}
\put(570, 0){\includegraphics[width=0.42\textwidth]{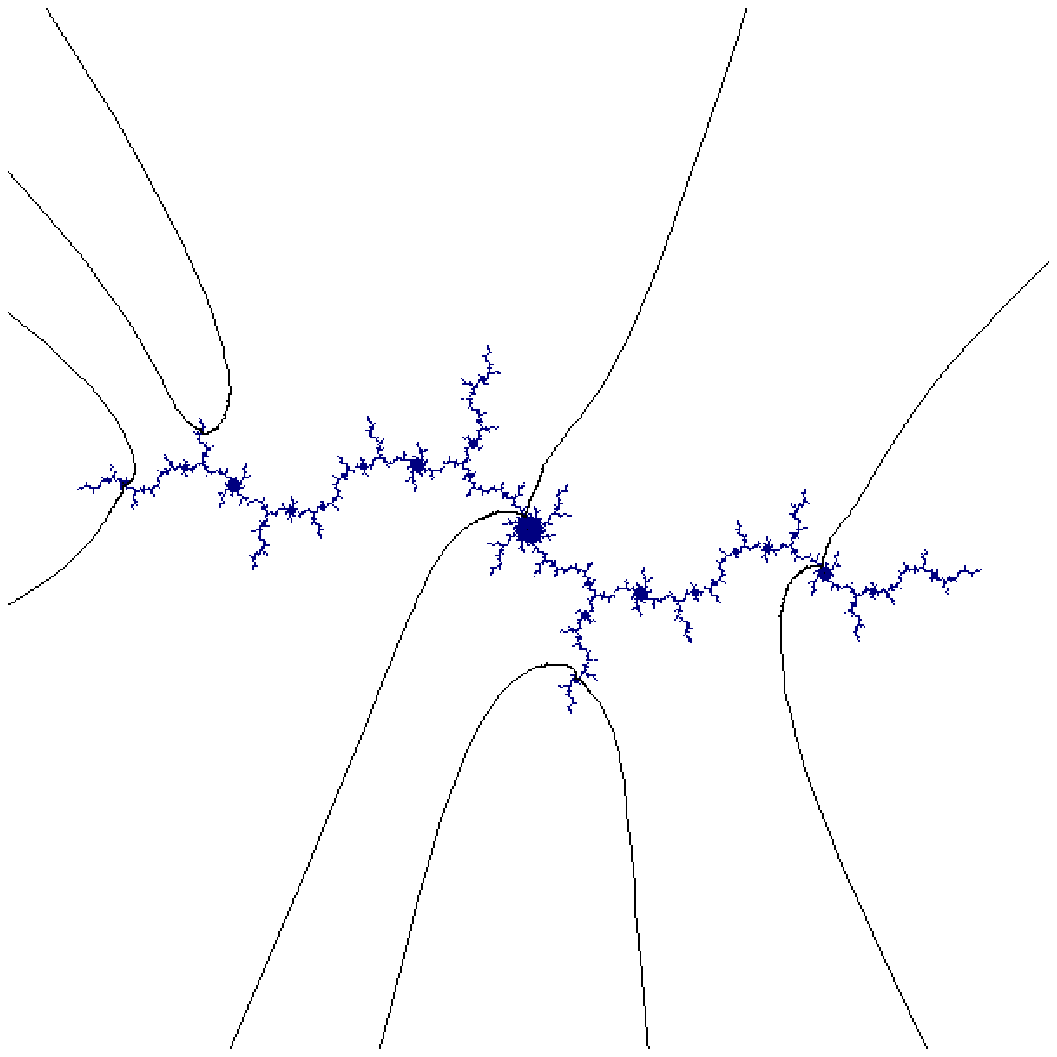}}
\thinlines
\multiput(10, 0)(560, 0){2}{\line(1, 0){420}}
\multiput(10, 420)(560, 0){2}{\line(1, 0){420}}
\multiput(10, 0)(560, 0){2}{\line(0, 1){420}}
\multiput(430, 0)(560, 0){2}{\line(0, 1){420}}
\put(445, 140){\line(3, 1){110}}
\put(445, 300){\line(1, 0){110}}
\put(445, 367){\line(6, -1){110}}
\end{picture} \caption[]{\label{Fchsh}
A formal mating $g=P\fmate Q$. The Julia set $\K_p$ for the five-periodic
center $p$ corresponding to $\gamma_\sM(11/31)$ is shown on the right; the
Misiurewicz Julia set $\K_q$ with $q=\gamma_\sM(19/62)$ in the left image is
rotated. (This does not change the set itself, but its external rays.) The
ray connections and dynamics are discussed in the main text.\\[1mm]
There are various ways to visualize the sets
$\phi_0(\K_p),\,\phi_\infty(\K_q)\subset\hat\C$ in the plane $\C$\,: instead of
$\K_q$ coming from $\infty$, we may rotate the sphere such that $\K_q$ is
translated above or below $\K_p$\,, or to save space here, translated to the
left or right and rotated.\\[1mm]
In any case, $\phi_0(\r_p(\theta))$ is connected with
$\phi_\infty(\r_q(-\theta))$; three connections are indicated between the two
images. When discussing the combinatorics of a ray-equivalence class, we may
avoid conjugation of several angles by assuming that $\r_p(\theta)$ connects to
$\r_{\overline q}(\theta)$, but to draw these rays without crossing, you would
need to stack two sheets of paper.}
\end{figure}

In general both $\K_p$ and $\K_q$ contain pinching points and branch points
with several rays landing together, so there are ray-equivalence classes
consisting of subsequent rays connecting points in $\partial\K_p$ and
$\partial\K_q$ alternately. For rational angles, the landing pattern is
understood combinatorially, and the identifications of periodic and preperiodic
points can be determined. Consider the example of the $5$-periodic $p$ with the
external angle $11/31$ and preperiodic $q$ with angle $19/62$ in
Figure~\ref{Fchsh}: since $q$ belongs to the $2/5$-limb of $\M$, there are five
branches of $\K_q$ and five external rays at the fixed point $\alpha_q$\,,
which are permuted with rotation number $2/5$ by $Q$. Now $p$ is chosen such
that the complex conjugate angles land pairwise with another $5$-cycle at the
Fatou basins; the rays of $Q$ corresponding to the latter angles land at
endpoints, including the iterate $Q(q)$ of the critical value. So in the
topological mating and in the geometric mating $f\eqg P\tmate Q$, the point
$Q(q)$ is identified both with $\alpha_q$ and with a repelling $5$-cycle of
$P$. Now the critical point $0$ of $f$ is $5$-periodic, while $f^2(\infty)$ is
fixed. The five components of the immediate attracting basin all
touch at this fixed point with rotation number $2/5$, although they had
disjoint closures in $\K_p$\,.

Basic definitions and the geometry of ray-equivalence classes are discussed
in \textbf{Section~\ref{2}}. Simple examples of shared matings and of mating
discontinuity are obtained in \textbf{Section~\ref{3}}. The rational map $f$
above belongs to the same one-parameter family as matings with the Chebychev
polynomial, but it is not of this form. There are five other representations as
a mating: take the Rabbit with rotation number $2/5$ for $P$ and suitable
preperiodic parameters $q_1\,,\,\dots,\,q_5$ for $Q$, which are related to
the angles at $-\alpha_p$\,. More generally, we have
$P\tmate Q_i=P\tmate Q_j$ for all $p$ in the small satellite Mandelbrot set,
since the rays at $-\alpha_p$ are homotopic with respect to the postcritical
set and so the precaptures are combinatorially equivalent. Taking higher
rotation numbers gives shared matings with larger multiplicity. While it is
obvious that a hyperbolic rational map has only a finite number of
representations as a mating, this is not known in general when one or both of
the critical points are preperiodic. Finiteness is shown here for Chebychev
maps with one critical point periodic, and in \cite{pmate} for
Latt\`es maps. Examples with arbitrarily high multiplicity are obtained as well
for matings of the Airplane with preperiodic polynomials; here preperiod and
period are of the same order as the multiplicity, in contrast to the hyperbolic
examples by Rees \cite{rap}, where the period grows exponentially. ---
Simple ray connections can be used to define preperiodic matings with
$f(0)=\infty$. This property is lost when preperiodic parameters
converge to a parabolic parameter, confirming that mating is not jointly
continuous. The mechanism is similar to geometrically infinite examples by
Bl\'e--Valdez--Epstein \cite{ble2, epstc}, but here all maps are geometrically
finite and matability does not require special arguments.

In general there is only a Cantor set of angles 
at the Hubbard tree $T_{\overline q}\subset\K_{\overline q}$\,, whose Hausdorff
dimension is less than $1$. If an open interval in the complement contains all
angles on one side of the arc $[-\alpha_p\,,\,\alpha_p]\subset\K_p$\,, ray
connections of the formal mating $P\fmate Q$ are bounded explicitly, and the
topological mating exists. This approach was used by Shishikura--Tan in a cubic
example \cite{st3}; in the quadratic case it generalizes the treatment of
$1/4\tmate1/4$ by Milnor \cite{m14} to large classes of examples. These include
the mating of Airplane and Kokopelli, answering a question by Adam Epstein
\cite{mateq}: can the mating be constructed without employing the 
theorems of Thurston and Rees--Shishikura--Tan? See \textbf{Section~\ref{4}}.
Note however, that only the branched covering on the glued Julia sets is
constructed here, not a conjugate rational map. On the other hand, the method
applies to geometrically infinite parameters as well. Examples
of irrational ray connections and an algorithm for finding long ray connections
are discussed in addition. In \textbf{Section~\ref{5}}, specific ray
connections for polynomials from conjugate limbs are obtained, which are
related to renormalization of one polynomial. These ray-equivalence classes
accumulate on the Julia set, such that the quotient space is not Hausdorff.

 --- This is the second paper in a series on matings and other applications
of the Thurston Algorithm \cite{teich, dh, book2h}:
\begin{itemize}
\item \textbf{The Thurston Algorithm for quadratic matings \cite{qmate}.}
The Thurston Algorithm for the formal mating is implemented by pulling back
a path in moduli space; an alternative initialization by a
repelling-preperiodic capture is discussed as well. When the Thurston
Algorithm diverges in ordinary Teichm\"uller space due to postcritical
identifications, it still converges on the level of rational maps and colliding
marked points --- it is not necessary to implement the essential mating by
encoding ray-equivalence classes numerically. The proof is based on the
extended pullback map on augmented Teichm\"uller space constructed by Selinger
\cite{ext, char}.
\item \textbf{Quadratic matings and ray connections} [the present paper].
\item \textbf{Quadratic matings and Latt\`es maps \cite{pmate}.}
Latt\`es maps of type \2 or \4 are represented by matings in basically nine,
respectively three, different ways. This is proved from combinatorics of
polynomials and ray-equivalence classes. The Shishikura Algorithm relates the
topology of the formal mating to the multiplier of the corresponding affine
map on a torus. The slow mating algorithm diverges in certain cases: while
the expected collisions are happening, a neutral eigenvalue from the
one-dimensional Thurston Algorithm persists, producing an attracting center
manifold in moduli space. (Joint work with Arnaud Ch\'eritat.) Twisted
Latt\`es maps are discussed as well, and the Hurwitz equivalence between
quadratic rational maps with the same ramification portrait is constructed
explicitly, complementing the approach related to the moduli space map by
Sarah Koch \cite{endo}.
\item \textbf{Slow mating and equipotential gluing \cite{emate}, jointly with
Arnaud Ch\'eritat.}
Equipotential gluing is an alternative definition of mating, not based on the
Thurston Algorithm. Equipotential lines of two polynomials are glued to
define maps between spheres, and the limit of potential $0$ is considered.
The initialization of the slow mating algorithm depends on an initial radius
$R$; when $R\to\infty$, slow mating is shown to approximate equipotential
gluing. The visualization in terms of holomorphically moving Julia sets and
their convergence is discussed and related to the notion of conformal
mating.
\item \textbf{Quadratic captures and anti-matings \cite{amate}.}
The slow Thurston Algorithm is implemented for captures and for anti-matings
as well. The latter means that two planes or half-spheres are mapped to each
other by quadratic polynomials, and the filled Julia sets of two quartic
polynomials are glued together.  There are results analogous to matings, but
a complete combinatorial description does not exists due to the complexity of
even quartic polynomials. For specific families of quadratic rational
maps, the loci of mating, anti-mating, and captures are obtained numerically.
\item \textbf{The Thurston Algorithm for quadratic polynomials \cite{smate}.}
The slow Thurston Algorithm is implemented for several kinds of Thurston maps
giving quadratic polynomials. These include a spider algorithm with a path
instead of legs, Dehn twisted polynomials, moving the critical value by
recapture or precapture, and tuning. Using the Selinger results on removable
obstructions, the spider algorithm is shown to converge in the obstructed
case of satellite Misiurewicz points as well. Recapture surgery is related
to internal addresses, and used to discuss a specific example of twisted
polynomials.
\end{itemize}
\textbf{Acknowledgment:}
Several colleagues have contributed to this work by inspiring discussions.
I wish to thank in particular
Laurent Bartholdi, Adam Epstein, Mikhail Hlushchanka, Daniel Meyer, Mary Rees,
Dierk Schleicher, and Tan Lei. And I am grateful to the mathematics department
of Warwick University for their hospitality.

\section{Mating: definitions and basic properties} \label{2}
After recalling basic properties of quadratic polynomials and matings, the
geometry of rational and irrational ray-equivalence classes is described,
generalizing an observation by Sharland \cite{shar2}. Repelling-preperiodic
captures are considered as an alternative construction of matings; the proof
was given in \cite{qmate}, using the relation between ray-equivalence classes
and Thurston obstructions from \cite{rst, rs}.

\subsection{Polynomial dynamics and combinatorics} \label{2p}
For a quadratic polynomial $f_c(z)=z^2+c$, the filled Julia set $\K_c$
contains all points $z$ with $f_c^n(z)\not\to\infty$. It is connected, if and
only if the critical point $z=0$ does not escape, and then the parameter $c$
belongs to the Mandelbrot set $\M$ by definition. A \textbf{dynamic ray}
$\r_c(\theta)$ is the preimage of a straight ray with angle $2\pi\theta$ under
the Boettcher conjugation
$\Phi_c:\hat\C\setminus\K_c\to\hat\C\setminus\overline\disk$. For rational
$\theta$, the rays and landing points are periodic or preperiodic under
$f_c$\,, since $f_c(\r_c(\theta))=\r_c(2\theta)$. If two or more periodic
rays land together, this defines a non-trivial orbit portrait; it
exists if and only if the parameter $c$ is at or behind a certain root
\cite{ser, mer}. There are analogous parameter rays with rational angles
$\r_\sM(\theta)$ landing at roots and Misiurewicz points; the angles of a root
are characteristic angles from the orbit portrait. In particular,
the $k/r$-limb and wake of the main cardioid are defined by two parameter
rays with $r$-periodic angles, and for the corresponding parameters $c$, the
fixed point $\alpha_c\in\K_c$ has $r$ branches and external angles permuted
with rotation number $k/r$. Denote landing points by
$z=\gamma_c(\theta)\in\partial\K_c$ and $c=\gamma_\sM(\theta)\in\partial\M$,
respectively. $f_c$ is \textbf{geometrically finite}, if it is preperiodic,
hyperbolic, or parabolic.

\begin{prop}[Douady Magic Formula, Bl\'e]\label{Pdmf}
Suppose $\theta\in[0,\,1/3]$ is an external angle of the main cardioid, then
$\Theta=1/2+\theta/4\in[1/2,\,7/12]$ is an external angle of the real axis
$\M\cap\R$.
\end{prop}

\textbf{Proof:} According to \cite{bsoos}, the orbit of $\theta$ under
doubling is confined to $[\theta/2\,,(1+\theta)/2]$. Now taking a suitable
preimage shows that the orbits of $\theta$ and $\Theta$ never enter
$((\theta+1)/4,\,(\theta+2)/4)\supset(1-\Theta,\,\Theta)$, so $\Theta$ is
combinatorially real: it defines a unique real parameter $c$ by approximation,
and the parameter ray $\r_\sM(\Theta)$ accumulates at a fiber \cite{sf2}
intersecting the real line in $c$. Bl\'e \cite{ble1} has shown that $f_c$ is
strongly recurrent but not renormalizable, so the fiber is trivial and the ray
actually lands, $c=\gamma_\sM(\Theta)$. \mybox

\subsection{Topological mating and geometric mating} \label{2d}
For parameters $p,\,q\in\M$ with locally connected Julia sets, define the
\textbf{formal mating} $g=P\fmate Q$ of the quadratic polynomials $P(z)=z^2+p$
and $Q(z)=z^2+q$ as follows: $g:\hat\C\to\hat\C$ is a branched covering with
critical points $0$ and $\infty$, and normalized such that $g(z)=z^2$ for
$|z|=1$. On the lower and upper half-spheres, $g$ is topologically conjugate
to $P$ and $Q$ by homeomorphisms $\phi_0$ and $\phi_\infty$\,, respectively.
An \textbf{external ray} $\r(\theta)$ of $g$ is the union of
$\phi_0(\r_p(\theta))$ and $\phi_\infty(\r_q(-\theta))$ plus a point on the
equator; each ray connects a point in $\phi_0(\K_p)$ to a point in
$\phi_\infty(\K_q)$. A \textbf{ray-equivalence class} is a maximal connected
set consisting of rays and landing points. Collapsing all classes to points may
define a Hausdorff space homeomorphic to the sphere; then the map corresponding
to $g$ is a branched covering again \cite{maten}, which defines the
\textbf{topological mating} $P\tmate Q$ up to conjugation. By the
identifications, periods may be reduced and different orbits meet. We are
interested in a rational map $f$ conjugate to the topological mating, and we
shall speak of ``the'' geometric mating when the following normalization is
used. 
Note however, that uniqueness is not obvious when the polynomials are
not geometrically finite, in particular if there is a locally connected
Julia set carrying an invariant line field.

\begin{dfn}[Normalization of the geometric mating]\label{Dgm}
Suppose the topological mating $P\tmate Q$ is topologically conjugate to a
quadratic rational map $F$, and the conjugation $\psi$ is conformal in the
interior of the filled Julia sets. Then the geometric mating exists
and it is M\"obius conjugate to $F$.

The \textbf{geometric mating} $f\eqg P\tmate Q$ is normalized such that $\psi$
maps the critical point of $P$ to $0$, the critical point of $Q$ to $\infty$,
and the common $\beta$-fixed point to $1$. If the latter condition is dropped,
then $f$ is affine conjugate to the geometric mating, and we shall write
$f\eqa P\tmate Q$.
\end{dfn}

Sometimes it is convenient to write $p\tmate q$ or $\theta_p\tmate\theta_q$
for $P\tmate Q$; here a periodic angle is understood to define a center, not
a root. In the postcritically finite case, the geometric mating is constructed
using Thurston theory as follows:

\begin{thm}[Rees--Shishikura--Tan]\label{Trstrs}
Suppose $P$ and $Q$ are postcritically finite quadratic polynomials, not from
conjugate limbs of the Mandelbrot set. Then the geometric mating
$f\eqg P\tmate Q$ exists.
\end{thm}

Idea of the \textbf{proof:}
The formal mating $g=P\fmate Q$ is a postcritically finite branched covering,
a Thurston map. So it is combinatorially equivalent to a rational map, if and
only if it is unobstructed, excluding type \2 here \cite{pmate}. According to
Rees--Shishikura--Tan, all obstructions are L\'evy cycles converging to
ray-equivalence classes under iterated pullback \cite{rst}. See the example in
Figure~3 of \cite{qmate}. In the case of non-conjugate limbs, these
obstructions are removed by collapsing postcritical ray-equivalence trees,
which defines an unobstructed essential mating $\tilde g$. Now the Thurston
Theorem \cite{teich, dh, book2h} produces a rational map $f$ equivalent to $g$
or $\tilde g$, respectively, unique up to normalization. By iterating a
suitable equivalence, a semi-conjugation from $g$ to $f$ is obtained \cite{rs},
which collapses all ray-equivalence classes to points. So $f$ is conjugate to
the topological mating $P\tmate Q$. \mybox 

\begin{conj}[Quadratic mating]\label{Cqm}
For quadratic polynomials $P$ and $Q$ with locally connected Julia sets, the
geometric mating exists, unless $p$ and $q$ are in conjugate limbs of the
Mandelbrot set.
\end{conj} 

Originally, it was expected that mating depends continuously on the
polynomials \cite{mquad}; 
various counterexamples by Adam Epstein \cite{epstc, mateq} are discussed in
Section~\ref{3c}, and a simple new counterexample is given. ---
The geometric mating is known to exist in the following quadratic cases:
\begin{itemize}
\item In the postcritically finite situation, Conjecture~\ref{Cqm} was proved
in \cite{rst, rs}, cf.~Theorem~\ref{Trstrs}. In this case, the geometric
mating exists, whenever the topological mating does.
See \cite{maten, emate} for various notions of conformal mating.
\item Suppose $P$ and $Q$ are hyperbolic quadratic polynomials, and denote the
corresponding centers by $p_0$ and $q_0$\,, let $f_0\eqg P_0\tmate Q_0$\,.
Now $P_0$ is quasiconformally conjugate to $P$ in a neighborhood of the
Julia set $\J_{p_0}=\partial\K_{p_0}$\,, analogously for $Q_0$\,, and there is
a rational map $f$ with the corresponding multipliers, such that $f_0$ is
quasiconformally conjugate to $f$ in a neighborhood of $\J_{f_0}$\,. The
conjugations of polynomials respect the landing of dynamic rays, so the
semi-conjugations from $P_0$ and $Q_0$ to $f_0$ define new semi-conjugations
from $P$ and $Q$ to $f$ in neighborhoods of the Julia sets. Using conformal
conjugations to Blaschke products on the immediate basins,
the required semi-conjugations from $\K_p\fmate\K_q\to\hat\C$ are constructed,
and $f\eqg P\tmate Q$ is a geometric mating. The same argument works when one
polynomial is hyperbolic and the other one is preperiodic.
\item A geometrically finite quadratic polynomial is preperiodic,
hyperbolic, or parabolic. Ha{\ddoti}ssinsky--Tan have constructed all matings
of geometrically finite polynomials from non-conjugate limbs \cite{htl}: when
parabolic parameters are approximated radially from within hyperbolic
components, the geometric matings converge. 
The proof is based on distortion control techniques by Cui. On the other hand,
when two parabolic parameters are approximated tangentially, mating may be
discontinuous; see \cite{epstc, mateq} and Section~\ref{3c}.
\item For quadratic polynomials having a fixed Siegel disk of bounded type,
Yampolsky--Zakeri \cite{yz} construct the geometric mating when the multipliers
are not conjugate, and obtain the mating of one Siegel polynomial with the
Chebychev polynomial in addition. The proof combines Blaschke product models,
complex a priori bounds, and puzzles with bubble rays.
\item Suppose $\theta$ defines a parameter $p$ with a Siegel disk of bounded
type and consider the real parameter $q$ with angle $\Theta=1/2+\theta/4$
defined in Proposition~\ref{Pdmf}, which is strongly recurrent. The geometric
mating $f\eqg P\tmate Q$ exists according to Bl\'e-Valdez \cite{ble2}.
\item Denote the family of quadratic rational maps $f_a(z)=(z^2+a)/(z^2-1)$
with a superattracting $2$-cycle by $V_2$\,. It looks like a mating between
the Mandelbrot set $\M$ and the Basilica Julia set $\K_B$\,, both truncated
between the rays with angles $\pm1/3$. Capture components correspond to
Fatou components of the Basilica. Large classes of maps in $V_2$ are known to
be matings of quadratic polynomials with the Basilica, by work of
Luo, Aspenberg--Yampollsky, Dudko, and Yang \cite{luo, aym, ddlam, yang}.
The basic idea is to construct puzzle-pieces with bubble rays both in the
dynamic plane and in the parameter plane. 
This approach does not seem to generalize to $V_3$\,, because Rabbit matings
may be represented by Airplane matings as well.
\item When $p$ is periodic and $\overline q$ shares an angle with a boundary
point of a preperiodic Fatou component, the geometric mating is constructed by
regluing a capture according to Mashanova--Timorin \cite{mtrg}.
\item For large classes of geometrically finite and infinite examples,
Theorem~\ref{Tbdrcxpl} shows that ray-equivalence classes
are uniformly bounded trees. So the topological mating exists according to
Epstein \cite{maten}, but the geometric mating is not constructed here.
\end{itemize}
In higher degrees, a topological mating $P\tmate Q$ may exist when there is no
geometric mating. An example with periodic cubic polynomials is discussed in
\cite{st3}. Other examples are obtained from expanding Latt\`es maps: choose a
$2\times2$ integer matrix $A$ with trace $t$ and determinant $d$ satisfying
$0<t-1<d<t^2/4$, e.g., $t=d=5$. This defines a Thurston map $g$
of type \2 with degree $d$. Now $g^n$ is expanding and not equivalent to a
rational map, since the eigenvalues of $A^n$ are real $>1$ and distinct
\cite{dh, book2h, pmate}. But according to \cite{meyerp}, $g^n$ is a
topological mating for large $n$. 

\subsection{Ray connections and ray-equivalence classes} \label{2r}
For the mating of quadratic polynomials $P(z)=z^2+p$ and $Q(z)=z^2+q$ with
locally connected Julia sets, rays and ray-equivalence classes are defined in
terms of the formal mating $g=P\fmate Q$\,.
A ray connection is an arc within a ray-equivalence class. The length of an
arc or loop is the number of rays involved, and the diameter of a
ray-equivalence class is the greatest distance with respect to this notion of
length. 
We shall discuss the structure of ray-equivalence classes in detail for
various examples, and show existence of the topological mating in certain
cases. By the Moore Theorem \cite{maten}, all ray-equivalence classes must
be trees and the ray-equivalence relation must be closed. For this the length
of ray connections will be more important than the number of rays and landing
points in a ray-equivalence class: there is no problem when, e.g., branch
points with an increasing number of branches converge to an endpoint,
since the angles will have the same limit. The following results are proved in
Propositions~4.3 and~4.12 of \cite{maten}:

\begin{prop}[Ray connections and matability, Epstein]\label{Palrcm}
Consider ray-equivalence classes for the formal mating $g=P\fmate Q$ of
$P(z)=z^2+p$ and $Q(z)=z^2+q$, with $\K_p$ and $\K_q$ locally connected.

$1$. If all classes are trees and uniformly bounded in diameter, the
topological mating $P\tmate Q$ exists as a branched covering of the sphere.

$2$. If there is an infinite or a cyclic ray connection, the topological
mating does not exist.
\end{prop}

Note that there is no statement about non-uniformly bounded trees. For
Misiurewicz matings having a pseudo-equator, Meyer \cite{meyerp} has
shown that ray-equivalence classes are bounded uniformly in size; hence the
diameters are bounded uniformly as well. Theorem~\ref{Tbdrcxpl} gives
topological matings $P\tmate Q$, where all ray-equivalence classes are
bounded uniformly in diameter, but they need not be bounded in size; see
Example~\ref{Xreqmb}. The following description of ray-equivalence classes
can be given in general, speaking of connections between $\partial\K_p$ and
$\partial\K_{\overline q}$ according to Figure~\ref{Fchsh}: 

\begin{prop}[Shape of ray-equivalence classes, following Sharland]\label{Psrec}
Consider rational and irrational ray-equivalence classes for the formal mating
$g=P\fmate Q$ of quadratic polynomials, with $\K_p$ and $\K_q$ locally
connected.

$1$. Any branch point of a ray-equivalence class is a branch point of $\K_p$
or $\K_{\overline q}$\,. Thus it is precritical, critical, preperiodic,
or periodic. So with countably many exceptions, all ray-equivalence classes are
simple arcs $($finite or infinite$)$, or simple loops.

$2$. Suppose the periodic ray-equivalence class $C$ is a finite tree, then all
the angles involved are rational of the same ray period $m$. Either $C$ is an
arc and $m$-periodic as a set, or it contains a unique point $z$ of period
$m'=m/r$ with $r\ge2$ branches. Then $z$ is the only possible branch point
of $C$, so $C$ is a topological star when $r\ge3$.

$3$. Suppose that the topological mating $P\tmate Q$ exists. Then only critical
and precritical ray-equivalence classes may have more than one branch point.
More precisely, we have the following cases:\\
$a)$ Both $P$ and $Q$ are geometrically finite. Then irrational ray-equivalence
classes of $g$ are finite arcs, and rational ray-equivalence classes may have
at most seven branch points.\\ 
$b)$ Precisely one of the two polynomials is geometrically finite. Then
irrational classes have at most one branch point, and rational classes may have
up to three.\\
$c)$ Both polynomials are geometrically infinite. Then irrational classes have
at most three branch points, and rational classes have at most one.
\end{prop}

Item~2 was used by Sharland \cite{shar1, shar2} to describe 
hyperbolic matings with cluster cycles. It is employed in Sections~4.3
and~6 of \cite{pmate} to classify matings with orbifold of
essential type \2, and here in Section~\ref{3t}.

\textbf{Proof:}
1. Since the rays themselves are not branched, the statement is immediate from
the No-wandering-triangles Theorem \cite{tgdr, ser} for branch points
of quadratic Julia sets.

2. Rational rays landing together have the same preperiod and ray period, and
only rational rays land at periodic and preperiodic points of a locally
connected Julia set. So they never land together with irrational rays.
Ray-equivalence classes are mapped homeomorphically or as a branched cover. If
a finite tree $C$ satisfies $g^{m'}(C)\cap C\neq\emptyset$ with minimal
$m'\ge1$, we
have $g^{m'}(C)=C$ in fact, and $C$ does not contain a critical point. Since
$g^{m'}$ is permuting the points and rays of $C$, there is a minimal $m\ge m'$,
such that $g^m$ is fixing all points and rays, and all angles are $m$-periodic.
Suppose first that $C$ contains a branch point $z$ with $r\ge3$ branches. It is
of satellite type, so its period is $m/r\ge m'$, and the $r$ branches are
permuted transitively by $g^{m/r}$. Thus all the other points are $m$-periodic,
and they cannot be branch points, because the first return map would not
permute their branches transitively. So $m'=m/r$. On the other hand, if $C$ is
an arc, then $g^{m'}$ is either orientation-preserving and $m=m'$, or
orientation-reversing and $m=2m'$. In the latter case, the number of rays must
be even, since each point is mapped to a point in the same Julia set, and the
point in the middle has period $m'=m/2$.

3) A periodic ray-equivalence class may contain a single branch point according
to item~2. In case a) a preperiodic class may contain two postcritical points
(from different
polynomials), and we have a pullback from critical value to critical point
twice. Each time the number of branch points may be doubled, and a new branch
point be created. This can happen only once in case b) and not at all in case
c). On the other hand, an irrational ray-equivalence class $C$ may contain only
critical and precritical branch points, and this can happen only when the
corresponding polynomial is geometrically infinite. Some image of $C$ contains
postcritical points instead of (pre-)critical ones, and it can contain only
one postcritical point from each polynomial, since it would be periodic
otherwise. So pulling it back to $C$ again gives at most three branch points.
Note that an irrational periodic class would be infinite or a loop,
contradicting the assumption of matability. \mybox

\subsection{Matings as repelling-preperiodic captures} \label{2c}
A Thurston map may be defined by shifting a critical value to a preperiodic
point along a path \cite{rees1}:

\begin{prop}[and definition]\label{Phrc}
Suppose $P$ is a postcritically finite quadratic polynomial and
$z_1\in\K_p$ is preperiodic and not postcritical. Let the new
postcritical set be $P_g=P_P\cup\{P^n(z_1)\,|\,n\ge0\}$. Consider an arc $C$
from $\infty$ to $z_1$ not meeting another point in $P_g$ and choose a
homeomorphism $\phi$ shifting $\infty$ to $z_1$ along $C$, which is the
identity outside off a sufficiently small neighborhood of $C$. Then:\\[1mm]
$\bullet$ $g=\phi\circ P$ is well-defined as a quadratic Thurston map with
postcritical set $P_g$\,. It is a \textbf{capture} if $z_1$ is eventually
attracting and a \textbf{precapture} in the repelling case.\\[1mm]
$\bullet$ The combinatorial equivalence class of $g$ depends only on the
homotopy class of the arc $C$.
\end{prop}

See also the discussion of a possible numerical implementation of the Thurston
Algorithm in \cite{qmate}. Motivated by remarks of Rees and Mashanova--Timorin
\cite{mtrg}, the following result provides an alternative construction of
quadratic matings in the non-hyperbolic case; see the proof of Theorem~6.3 in
\cite{qmate}:

\begin{thm}[Matings as precaptures, following Rees]\label{TRprec}
Suppose $P$ is postcritically finite and $\theta$ is preperiodic, such that
$q=\gamma_\sM(-\theta)$ is not in the conjugate limb and
$z_1=\gamma_p(\theta)\in\partial\K_p$ is not postcritical. Then the precapture
$g_\theta=\phi_\theta\circ P$ along $\r_p(\theta)$ is combinatorially
equivalent or essentially equivalent to the geometric mating $f$ defined by
$P\tmate Q$.
\end{thm}

\section{Mating as a map between parameter spaces} \label{3}
Mating provides a partial map from $\M\times\M$ to the moduli space of
quadratic rational maps. This map is neither surjective, injective,
nor continuous. The characterization of matings in terms of equators and
pseudo-equators by Thurston--Wittner and Meyer is discussed in
Section~\ref{3e}. Old and new examples of shared matings are described in
Section~\ref{3s}, and particular sequences with arbitrarily high
multiplicity are obtained in Sections~\ref{3t} and~\ref{3u}. Epstein has
given various examples of mating discontinuity, which are described
in Section~\ref{3c}, and a simple new construction is presented.

\subsection{Characterization of matings} \label{3e}
Hyperbolic quadratic rational maps $f$ are classified as follows according to
Rees \cite{rees1} and Milnor \cite{mquad}:\\[1mm]
B or II is \textbf{bitransitive}: both critical points are in the same cycle
of Fatou components but not in the same component.\\[1mm]
C or III is a \textbf{capture}: one critical point is in a strictly
preperiodic Fatou component.\\[1mm]
D or IV has \textbf{disjoint} cycles of Fatou components.\\[1mm]
E or I is \textbf{escaping}: both critical orbits converge to a fixed point
within the only Fatou component.

Now each hyperbolic component of type B, C, D contains a unique postcritically
finite map up to normalization, but there is no such map of type E.
While hyperbolic anti-matings may be of type B, C, or D \cite{amate}, every
hyperbolic mating is of type D. The converse is false according to
Ben Wittner \cite{matew}: 

\begin{xmp}[Wittner]\label{Xw34}
There is a unique real quadratic rational map of the form
$f_w(z)=(z^2+a)/(z^2+b)$\,, such that $0$ is four-periodic and $\infty$ is
three-periodic; approximately $a=-1.3812$ and $b=-0.3881$. 
This map is not a geometric mating of quadratic polynomials.
\end{xmp} 

\textbf{Proof:} Any mating $f\eqa P\tmate Q$ has this branch portrait, if
and only if $P$ is four-periodic and $Q$ is three-periodic. Wittner determined
all combinations numerically and found them to be different from $f_w$\,.
Alternatively, show combinatorially that
all of these matings have periodic Fatou components with common boundary
points; this is obvious when $P$ or $Q$ is of satellite type. Otherwise
$P$ or $\overline P$ is the Kokopelli at $\gamma_\sM(3/15)$ and $Q$ is the
Airplane --- then four-periodic Fatou components are drawn together pairwise
by ray-connections through the two-cycle of the Airplane. On the other hand,
for $f_w$ no closed Fatou components in the same cycle meet, since the critical
orbits are ordered cyclically as $z_3<w_2<z_0<z_2<w_1<z_1<w_0$ on
$\R\cup\{\infty\}$. 
The Julia set $\J_w$ is a Sierpinski carpet in fact \cite{mquad}. \mybox

The characterization of matings by an equator is a folk theorem going back to
Thurston; it was proved in \cite{matew, meyeru} under similar assumptions.
Statement and proof require some standard notions from Thurston
theory, see \cite{teich, dh, book2h, qmate}.

\begin{thm}[Thurston--L\'evy--Wittner]\label{Tmeqw}
Suppose $f$ is a postcritically finite rational map of degree $d\ge2$. Then
$f$ is combinatorially equivalent to a formal mating $g=P\fmate Q$, if and
only if it has an \textbf{equator} $\gamma$: a simple closed curve with the
property that $\gamma'=f^{-1}(\gamma)$ is connected and homotopic to
$\gamma$ relative to the postcritical set, traversed 
in the same direction. 
\end{thm}

\textbf{Proof:} By construction, a formal mating $g=P\fmate Q$ has the
equator $S^1$. So if $f$ is combinatorially equivalent to $g$, with 
$\psi_0\circ g=f\circ\psi_1$\,, then $\gamma=\psi_0(S^1)$ is homotopic to
$\gamma'=f^{-1}(\gamma)=\psi_1(S^1)$. Conversely, when $f$ has an equator, it
is equivalent to a Thurston map $\hat g$ with $\hat g(z)=z^d$ for $z\in S^1$.
So $\hat g$ is a formal mating of two topological polynomials $\hat P$ and
$\hat Q$. Suppose $\hat P$ is obstructed, thus $f$ is obstructed as well, then
it would be a flexible Latt\`es map with four postcritical points. Now $\hat P$
and $\hat Q$ together have six postcritical points; since $\hat P$ has at least
four, $\hat Q$ has at most two, so $\hat Q$ and $f$ have a periodic critical
point. But Latt\`es maps have only preperiodic critical points, so $\hat P$ and
$\hat Q$ are unobstructed in any case. By the Thurston Theorem, there are
equivalent polynomials $P$ and $Q$, which are determined uniquely by requiring
them monic, centered, and with suitable asymptotics of the $0$-ray under the
equivalence. 
Now $\hat g$ is equivalent to the formal mating $g=P\fmate Q$. \mybox

\begin{rmk}[Equator and pseudo-equator]\label{Rpeq}
1. Suppose $f\eqg P\tmate Q$ is a postcritically finite geometric mating. If
$f$ is hyperbolic, it is combinatorially equivalent to the formal mating
$g=P\fmate Q$, so it has an equator. If $f$ is not hyperbolic, there may be
identifications from postcritical ray-equivalence classes, such that $g$ is
obstructed and $f$ is combinatorially equivalent to an essential mating
$\tilde g$. Then $f$ does not have an equator corresponding to this
representation as a mating.

2. When $P$ and $Q$ have only preperiodic critical points, the essential
mating $\tilde g$ and the geometric mating $f$ may have a
\textbf{pseudo-equator}, which passes through all postcritical points; see
\cite{meyerp, meyeru} for the definition. The equator of $g$ is deformed to a
pseudo-equator of $\tilde g$, if and only if there are at most direct
ray connections between postcritical points. Conversely, when $f$ has a
pseudo-equator $\gamma$, each pseudo-isotopy from $\gamma$ to $f^{-1}(\gamma)$
determines a pair of polynomials $P,\,Q$ with $f\eqa P\tmate Q$.

3. A Thurston map $g$  is \textbf{expanding}, if there is a curve $C$ through
the postcritical points, such that its $n$-th preimages form a mesh with
maximal diameters going to $0$.
See \cite{bookbm, bookhp, bd4} for other notions of expansion. 
According to \cite{meyerp}, some iterate $g^n$ has a pseudo-equator and it
is equivalent to a topological mating. A finite subdivision rule may be used
to define an expanding map \cite{exp}; for an essential mating with a
pseudo-equator, Wilkerson \cite{mw1} constructs a subdivision rule from the
Hubbard trees.
\end{rmk}

\subsection{Shared matings} \label{3s}
A shared mating is a geometric mating with different representations,
$P_1\tmate Q_1\eqa f\eqa P_2\tmate Q_2$ with $P_1\neq P_2$ or $Q_1\neq Q_2$\,.
There are the following examples of shared matings, and techniques for
constructing them:
\begin{itemize}
\item Wittner \cite{matew} introduced the notion of shared matings and
discussed them for $V_3$ in particular. A simple example is given by the
geometric mating of Airplane and Rabbit, which is affine
conjugate to the geometric mating of Rabbit and Airplane,
$A\tmate R\eqa R\tmate A$. (Moreover, it is conjugate to a symmetric map,
which is not a self-mating.) Since the two polynomials are interchanged, this
example is called the Wittner flip. It can be explained by finding two
different equators, which has a few generalizations:
\item Exall \cite{exall} constructs pairs of polynomials $P,\,Q$ with
$P\tmate R\eqa Q\tmate A$ from a second equator. Using symbolic dynamics,
this can be done algorithmically. 
\item Rees \cite{rap} uses symbolic dynamics again to obtain unboundedly shared
Airplane matings. The period grows exponentially with the multiplicity.
\item Denote the rabbit of rotation number $k/n$ by $R$. There are $n-2$
primitive hyperbolic polynomials $Q$ of period $n$, such that $\overline Q$
has a characteristic angle from the cycle of $\alpha_r$\,. Then the rational
map $f\eqg R\tmate Q$ has a cluster cycle: both $n$-cycles of Fatou components
have a common boundary point, which is a fixed point corresponding to
$\alpha_r$\,. Tom Sharland \cite{shar1, shar2} has shown that $f$ is
determined uniquely by the rotation number and the relative displacement of
the critical orbits; $f$ has precisely two representations as a mating, which
are of the form $f\eqg R\tmate Q\eqa P\tmate R$.
\end{itemize}
When $f$ is a Latt\`es map, different representations are known
except in the case c) of $1/6\tmate1/6$. The Shishikura Algorithm can be
used to identify the particular map $f$ in the case of type \2,
and we have only one quadratic map of type \4. Combinatorial
arguments show that there are basically nine, respectively three, matings of
these types; see Sections~4 and~6 in \cite{pmate}.
\begin{itemize}
\item Case a) of type \2 is $1/4\tmate1/4 \eqa 23/28\tmate13/28
 \eqa 13/28\tmate23/28 \eqa 53/60\tmate29/60 \eqa 29/60\tmate53/60$.
\item Case b) of type \2 is given by
 $1/12\tmate5/12 \eqa -1/12\tmate5/12$.
\item Case d) of type \2 is
 $1/6\tmate5/14 \eqa 5/14\tmate 1/6 \eqa 3/14\tmate3/14 \eqa 3/14\tmate1/2
 \eqa 1/2\tmate3/14 \eqa 5/6\tmate1/2 \eqa 1/2\tmate5/6$.
\item Type \4 is given by
 $\pm1/4\tmate1/2 \eqa 5/12\tmate\pm1/6 \eqa 13/28\tmate\pm3/14$.
\end{itemize}
The following technique for producing shared matings is based on the
representation of matings as repelling-preperiodic captures according to
Theorem~\ref{TRprec}.

\begin{prop}[Shared matings from precaptures]\label{Pshprec}
Suppose $P(z)=z^2+p$ is geometrically finite, with $p\neq-2$, $p\neq 1/4$,
and $p$ not in the main cardioid. There are
countably many pairs of preperiodic angles $\theta_1\,,\,\theta_2$ such that:
the corresponding dynamic rays land together at a preperiodic pinching
point $z_1\in\partial\K_p$\,, which is not postcritical and not in the same
branch at $\alpha_p$ as $p$, and the branch or branches of $z_1$ between these
rays do not contain postcritical points of $P$ or iterates of $z_1$. Then we
have $P\tmate Q_1\eqa P\tmate Q_2$ with $q_i=\gamma_\sM(-\theta_i)$.
Moreover, $P\tmate Q_1\eqg P\tmate Q_2$ if $\beta_p$ is not between these rays.
\end{prop}

\textbf{Proof:} We need to exclude $p=1/4$ and the main cardioid, because
$\K_p$ would have no pinching points, and $p=-2$, because rays landing together
at the interval $\K_{-2}$ are never homotopic with respect to the postcritical
set. If $P$ is postcritically finite, Proposition~\ref{Phrc} shows that
the precaptures $\phi_{\theta_1}\circ P$ and $\phi_{\theta_2}\circ P$ are
combinatorially equivalent. So the canonical obstructions and the essential
maps are equivalent as well. According to the proof of Theorem~\ref{TRprec},
given in \cite{qmate}, the essential maps are equivalent to the geometric
matings. 
By continuity according to Section~\ref{2d}, the result extends to
geometrically finite $P$: \mybox
\begin{itemize}
\item The example $11/24\tmate13/56\eqg11/24\tmate15/56$ enjoys the following
property: the latter mating has an equator and a simple pseudo-equator,
while the former does not have either.
\item As another example, consider $p=\gamma_\sM(59/240)$ and
$q=\gamma_\sM(63/240)$. Applying this construction to $P$ and to $Q$ gives
$P\tmate P\eqg P\tmate Q$ and $Q\tmate P\eqg Q\tmate Q$. Here the first and
second polynomials may be interchanged on both sides, so we have four
representations of the same rational map; in particular there are shared
self-matings $P\tmate P\eqg Q\tmate Q$, and the flipped matings
$P\tmate Q\eqg Q\tmate P$.
\item When $P$ is the Basilica, all pinching points are preimages of
$\alpha_p$\,. Since none of these is iterated behind itself, shared matings are
obtained from any pinching point $z_1$\,, which is not $\alpha_p$ or behind it.
Dudko \cite{ddlam} has shown that these are the only shared Basilica matings,
since the parameter space is described as a mating of $\M$ and $\K_p$\,.
The simplest example is given by $P\tmate(z^2\pm\i)$: the geometric matings
are distinct and complex conjugate, and both affine conjugate to
$\frac{z^2+2}{z^2-1}$. The example $P\tmate5/24\eqa P\tmate7/24$ is
illustrated with a video of slow mating on
\href{http://www.mndynamics.com/index2.html}{www.mndynamics.com}\ .
Aspenberg \cite{aspbas} constructs the semi-conjugation from the Basilica
to the rational map, beginning with the Boettcher map;
in this alternative approach, shared matings are obtained from a non-unique
labeling of Fatou components by bubble rays.
\item Shared matings in the family of Chebychev maps are discussed in
Section~\ref{3t}. In certain cases, lower bounds on the multiplicity are
obtained from homotopic rays according to Proposition~\ref{Pshprec}, or
upper bounds are obtained directly.
\item When $z_1$ is a branch point of $\K_p$\,, there may be more than two
parameters $q_i$\,. In Theorem~\ref{Tunbampp} of Section~\ref{3u},
unboundedly shared Airplane matings with small preperiods and periods are
constructed. Although the Airplane does not contain any branch point, this is
achieved by choosing $q_i$ with a common branch point in $\K_q$\,.
\item If $f$ is a critically preperiodic rational map of degree $d\ge2$ with
three or four postcritical points, a pseudo-equator may produce several
unmatings by choosing different pseudo-isotopies to its preimage
\cite{meyeru}. A higher multiplicity is obtained when there are
degenerate critical points, or when a critical point is mapped to another one.
Probably the only quadratic example is the Latt\`es map of type \4.
See \cite{fkklpps} for related results on NET maps.
\end{itemize}
\begin{rmk}[Finite multiplicity]\label{Rfmsm}
If $f$ is a postcritically finite quadratic rational map, can there be
infinitely many representation as a mating $f\eqa P\tmate Q$?\\
$\bullet$ When $f$ is hyperbolic, there are only finitely many candidates
for $P$ and $Q$, since there are only finitely many quadratic polynomials with
a given superattracting period.\\
$\bullet$ When one critical point is periodic and one is
preperiodic, finiteness is not obvious. For a specific family, finiteness is
shown in Theorem~\ref{Tchmun} of the following section, using similar
techniques as in the Latt\`es case.\\ 
$\bullet$ When both critical points are preperiodic, finiteness is
shown for Latt\`es maps in \cite{pmate}. Probably the techniques can be applied
to a few other examples of small preperiod and period, but a general proof
shall be harder.
\end{rmk}

\subsection{Shared matings in the Chebychev family} \label{3t}
Let us define a Chebychev map as a quadratic rational map of the form
$f(z)=f_a(z)=\frac{z^2-a-2}{z^2+a}$\,, $a\neq-1$, for which $f(\infty)$ is
pre-fixed: $\infty\Rightarrow1\to-1\uparrow$. This family contains matings
with the Chebychev polynomial in particular:

\begin{prop}[Chebychev maps as matings]\label{Pchmch}
Suppose $P(z)=z^2+p$ and $Q(z)=z^2+q$ are geometrically finite and not in
conjugate limbs of the Mandelbrot set $\M$. Then the geometric mating is affine
conjugate to a Chebychev map, $f_a\eqa P\tmate Q$, if and only if $P$ and $Q$
are one of the following forms:\\[1mm]
a$)$ $Q$ is the Chebychev polynomial $Q(z)=z^2-2$ and $p$ is not in the
$1/2$-limb of $\M$.\\[1mm]
b$)$ $p$ is in the $k/r$-limb of $\M$, and $q=\gamma_\sM(-\theta)$, where
$\theta$ is one of the $r$ angles at $-\alpha_p\in\K_p\ ($which depend
only on $k/r)$.\\[1mm]
c$)$ For a rotation number $k/r\neq1/2$, denote the angles of the $k/r$-wake
by $\theta_\pm$ and let $\theta=(\theta_-+\theta_+)/2$ be the unique angle of
preperiod $1$ and period $r$ in that limb. If $q=\gamma_\sM(\theta)$, then $P$
must be in the closed wake of the primitive hyperbolic component $\Omega$ with
the root $\gamma_\sM(-2\theta)$.
\end{prop}

The Petersen transformation \cite{mquad} maps symmetric rational maps to
Chebychev maps, such that self-matings are mapped to Chebychev matings; see
also Remark~4.4 in \cite{pmate}. In the previous section the example of
shared self-matings $59/240\tmate59/240\eqg63/240\tmate63/240$ was obtained
from Proposition~\ref{Pshprec}; now the Petersen transformation gives the
shared Chebychev mating $59/240\tmate1/2\eqg63/240\tmate1/2$.

\textbf{Proof of Proposition~\ref{Pchmch}:} As explained in
Figure~\ref{Fchsh}, instead of saying that angles of $z\in\K_p$ and
$w\in\K_q$ are complex conjugate, we may say that $z\in\K_p$ shares an angle
with $\overline w\in\K_{\overline q}$\,, or connect $\K_{\overline p}$ to
$\K_q$ as well. In the formal mating $g=P\fmate Q$, the ray-equivalence class
of $g^2(\infty)$, corresponding to $\overline Q^2(0)=\overline Q(\overline q)$,
is fixed. By Proposition~\ref{Psrec}, it must contain a fixed
point of $P$ or $\overline Q$. If this is $\beta_p$ or $\beta_{\overline q}$\,,
the fixed class is the $0$-ray and
$\overline Q(\overline q)=\beta_{\overline q}$\,, which is case a).

b) Now suppose that $\overline Q(\overline q)$ is in the same ray-equivalence
class as $\alpha_p$ and $p\in\M_{k/r}$\,. Then the critical value $\overline q$
is connected to $-\alpha_p$\,. This connection must be direct, since $\K_p$
does not contain another pinching cycle of ray period $r$. So $\overline q$
shares an external angle with $-\alpha_p$ and all of these angles may occur,
since none is in the same sector at $\alpha_p$ as the critical value $p$, and
$q$ is not in the conjugate limb. The $r$ angles belong to different
Misiurewicz points in fact, since otherwise some $P\fmate Q$ would have a
closed ray connection.

c) Consider $q\in\M_{k/r}$ and $\overline P$ such that $Q(q)$ is in the same
ray-equivalence class as $\alpha_q$\,. The points are not equal, because the
preperiod would have to be $\ge r>1$. So the ray connection must have length
$2$, since length $\ge4$ would require additional pinching cycles of ray period
$r$ in $\M_{k/r}$\,. Thus $q$ has the external angle $\theta$ defined above,
and $\K_{\overline p}$ must contain a pinching cycle of period and ray period
$r$, which connects the cycle of $2\theta=\theta_-+\theta_+$ to that of
$\theta_\pm$\,. This cycle of $\K_{\overline p}$ persists from a primitive
hyperbolic component $\overline\Omega$ before $\overline p$.

It remains to show that $\overline\Omega$ exists and is unique. In the dynamic
plane of $Q$, the $r$ rays landing at $\alpha_q$ define $r$ sectors
$W_1\,,\dots,\,W_r$ with $0\in W_r$ and $q\in W_1$\,, such that $Q$ is a
conformal map $W_1\to W_2\to\dots\to W_{r-1}\to W_r$ and the sectors are
permuted with rotation number $k/r$. The external rays with angles
$2^{i-1}\theta_\pm$ bound $W_i$ for $1\le i\le r$. Now $W_i$ contains
$2^{i-1}\theta$ as well for $2\le i\le r-1$ and $W_r$ has both $2^{r-1}\theta$
and $2^r\theta=-\theta$. For $r\ge3$ it follows that $2\theta$ has exact period
$r$. We are looking for a primitive orbit portrait \cite{mer} 
connecting each angle in $\{2^i\theta\,|\,1\le i\le r\}$ to a unique angle in
$\{2^i\theta_-\,|\,1\le i\le r\}=\{2^i\theta_+\,|\,1\le i\le r\}$.

Starting in $W_r$\,, connect $2^{r-1}\theta$ to either $2^{r-1}\theta_-$ or
to $2^{r-1}\theta_+$\,, such that $2^r\theta$ is not separated from the other
angles. Pull the connection back until $2\theta$ is connected to $2\theta_-$
or $2\theta_+$\,. The complement of the $r-1$ disjoint small sectors is
connected, so we can connect the remaining angles $2^r\theta$ and $\theta_+$
or $\theta_-$ as well. This construction gives a valid orbit portrait and
defines $\overline\Omega$, which has the external angles $2\theta$ and
$2\theta_-$ or $2\theta_+$\,. Note that it is a narrow component, i.e., its
angular width is $1/(2^r-1)$ and there is no component of period $\le r$ behind
$\overline\Omega$. To show that $\overline\Omega$ is unique, suppose we
had started by connecting $2^{r-1}\theta$ with an angle not bounding $W_r$ and
pulled it back. This pullback would follow the rotation number $k/r$ as well
and the small sectors would overlap, the leaves would be linked. \mybox

Case b) provides maps from limbs of $\M$ to the Chebychev family, which are
partially shared according to Proposition~\ref{Pshprec}: e.g., for
$P$ geometrically finite in the $1/2$-limb, consider the geometric matings
corresponding to $P\tmate\pm1/6$, i.e. $p\mapsto f_a\eqa P\tmate\pm1/6$.
These two maps agree on the small Mandelbrot set of period $2$, but in
general do not agree on its decorations. Likewise, for $p$ in the $1/3$-limb,
we have three maps corresponding to $P\tmate3/14$,\, $P\tmate5/14$,\, and
$P\tmate13/14$,\, which agree on the small Mandelbrot set of period $3$. In
the decorations, two of the maps may agree on certain veins, but in general
the third one will be different: the relevant rays are no longer homotopic.
Note that according to case c), some of these Chebychev maps are represented
by $\tilde P\tmate3/14$ as well, with $\tilde p$ in the Airplane wake. In
particular, we have
$1/7\tmate3/14\eqg1/7\tmate5/14\eqa1/7\tmate13/14\eqa3/7\tmate3/14$.
Under the Petersen transformation mentioned above, this Chebychev map is the
image of $1/7\tmate3/7\eqa3/7\tmate1/7$, which is a symmetric map but not a
self-mating.

\begin{thm}[Chebychev maps as shared matings]\label{Tchmun}
Matings $P\tmate Q$ in the Chebychev family with hyperbolic $P$ have
non-uniformly bounded multiplicity:

$1$. Suppose $f=f_a$ is a Chebychev map, such that $z=0$ is $n$-periodic. Then
there are at most a finite number of representations $f_a\eqa P\tmate Q$.

$2$. For each rotation number $k/r$, there is a unique Chebychev map
$f=f_a$\,, such that $z=0$ is $r$-periodic and the fixed point $-1=f^2(\infty)$
is a common boundary point of the $r$ immediate basins, which are permuted with
rotation number $k/r$. This map has precisely $r+1$ realizations as a geometric
mating, $f\eqa P\tmate Q$, when $r\ge3$; for $r=2$ there are only $2$
representations. 
\end{thm}

\textbf{Proof:} 1. $P$ will be $n$-periodic, so there are only finitely many
possibilities for $P$. We must see that $r$ is bounded in cases b) and c).
But in both cases we have $r\le n$, since the wakes of period $r$ are narrow:
in case b) this is a basic property of limbs, and in case c) it was noted in
the proof of Proposition~\ref{Pchmch}.

2. In case a) of Proposition~\ref{Pchmch}, $z=-1$ corresponds to the
ray-equivalence class of angle $0$, which does not touch a hyperbolic
component of $P$. In cases b) and c), the rotation number at $-1$ is
precisely $k/r$, so the value of $k/r$ from the proposition must be the same
as in the hypothesis of the theorem; case c) is excluded for $k/r=1/2$.
In both cases, there is only one hyperbolic component of period $r$ in the limb
or wake. It remains to show that $f_a$ is unique, so that the $r+1$ (or two)
matings actually give the same map. Intuitively, this
follows from the fact that the hyperbolic component of $f_a$ bifurcates from
the hyperbolic component where $-1$ is attracting; the multiplier map with
$\rho=\frac{-4}{a+1}$ is injective for $|a+1|\ge4$. It can be proved by
Thurston rigidity, since there is a forward-invariant graph connecting the
postcritical points, which depends only on $k/r$ up to isomorphy. So all
possible maps $f_a$ are combinatorially equivalent, affine conjugate,
and equal. --- Note that the case of $k/r=2/5$ was discussed in the
Introduction and in Figure~\ref{Fchsh}. \mybox

\subsection{Unboundedly shared Airplane matings} \label{3u}
Denoting the Rabbit by $R$ and the Airplane by $A$, we have seen in the
previous Section~\ref{3t} that
$R\tmate3/14\eqg R\tmate5/14\eqa R\tmate13/14\eqa A\tmate3/14$. This example
belongs both to the Chebychev family and to the family $V_3$ with a
$3$-periodic critical point. Unboundedly shared matings were obtained in
Theorem~\ref{Tchmun}.2 by increasing both the period of the hyperbolic
polynomial $P$ and the ray period of the Misiurewicz polynomial $Q$. Another
example is obtained below, where $Q$ is always the Airplane, and the preperiod
of $P$ is unbounded. The proof will be a simple application of
Proposition~\ref{Pshprec} again. Airplane matings with unbounded multiplicity
are due to Rees \cite{rap} with hyperbolic polynomials $P$, such that the
period of $P$ grows exponentially with the multiplicity. 

\begin{thm}[Unboundedly shared Airplane matings]\label{Tunbampp}
For the Airplane $q$ and $n=3,\,5,\,7,\,\dots$, there are $n$ Misiurewicz
parameters $p_*\,,\,p_2\,,\,\dots,\,p_n$\, such that the geometric matings
agree, $f\eqg P_i\tmate Q$ for all $i=*,\,2,\,\dots,\,n$. Here all $p_i$ have
preperiod $n+1$, $p_*$ has period $1$ and $p_2\,,\,\dots,\,p_n$ have period
$n$; so $f(\infty)$ has preperiod $n+1$ and period $1$. The statement remains
true for large $n$, when $q$ is any geometrically finite parameter behind
$\gamma_\sM(5/12)$ and before the Airplane. E.g., $q$ may be the Misiurewicz
point $\gamma_\sM(41/96)$ as well.
\end{thm}

\begin{figure}[h!t!b!]
\unitlength 0.001\textwidth 
\begin{picture}(990, 420)
\put(10, 0){\includegraphics[width=0.42\textwidth]{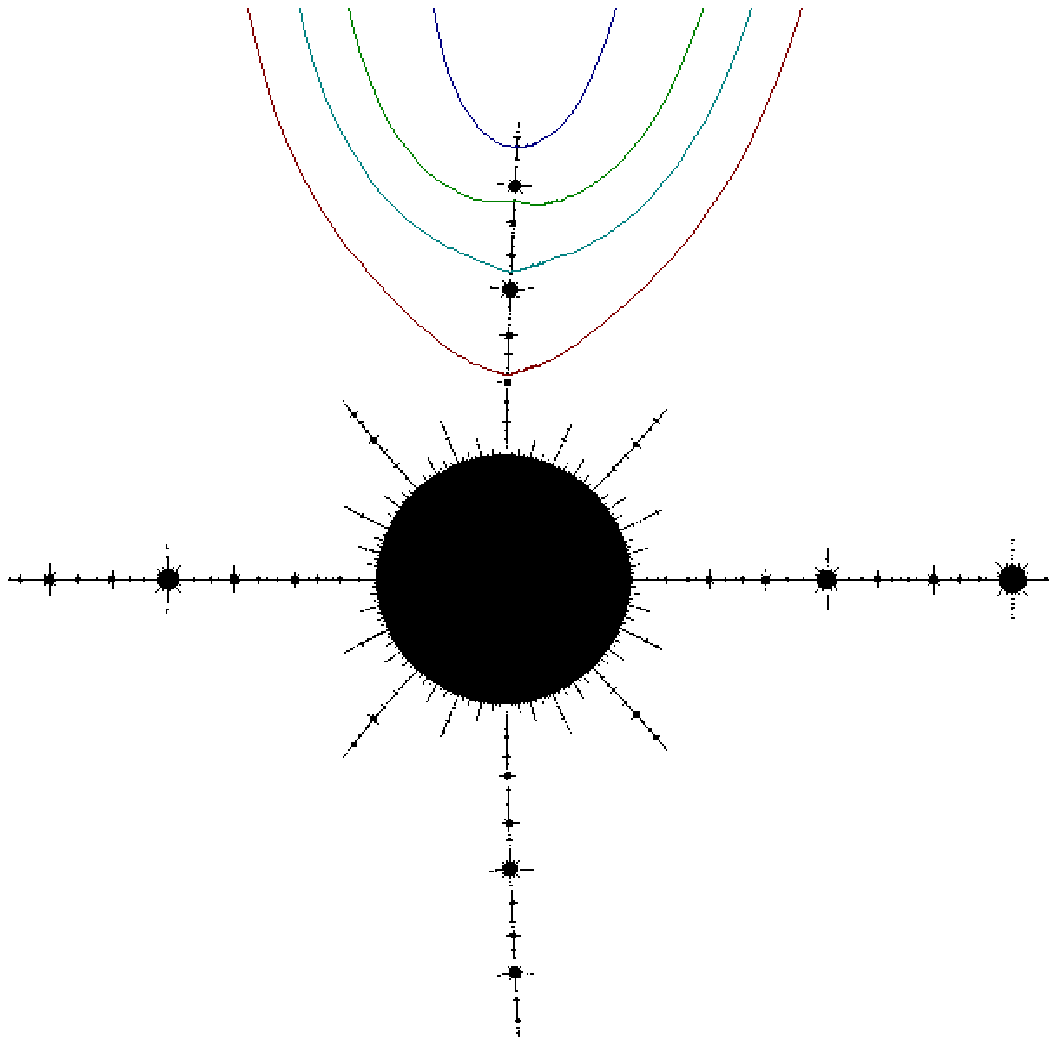}}
\put(570, 0){\includegraphics[width=0.42\textwidth]{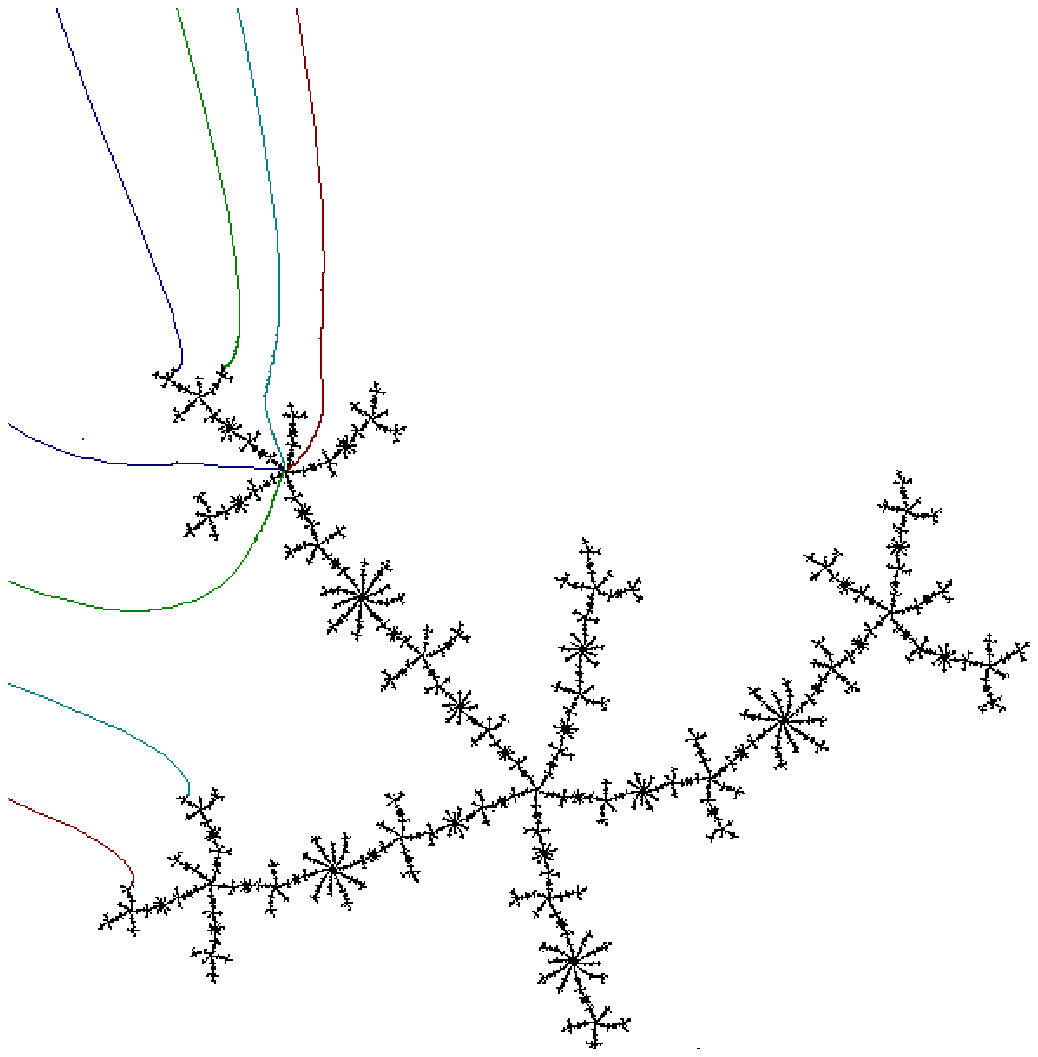}}
\thinlines
\multiput(10, 0)(560, 0){2}{\line(1, 0){420}}
\multiput(10, 420)(560, 0){2}{\line(1, 0){420}}
\multiput(10, 0)(560, 0){2}{\line(0, 1){420}}
\multiput(430, 0)(560, 0){2}{\line(0, 1){420}}
\end{picture} \caption[]{\label{Fmult}
Consider the formal mating $g=P\fmate Q$, with the Airplane $\K_q$ shown
rotated on the left, and $\K_p$ on the right for some $p$ in the $2/5$-limb.
According to the proof of Theorem~\ref{Tunbampp}, there are eight angles
$\theta_2\,,\,\dots,\,\theta_5\,,\,\theta_2'\,,\,\dots,\,\theta_5'$\,, such
that $-\theta_i$ and $-\theta_i'$ land together at the Airplane
$\partial\K_q$\,, while $\theta_i$ land together at $\partial\K_p$\,. So the
eight rays belong to a preperiodic ray-equivalence class of diameter four;
actually there are two more rays crossing the Airplane on the real axis.
Now there are five parameters $p=p_*\,,\,p_2\,,\,\dots,\,p_5$\,, such that
this ray-equivalence class contains the critical value $p$, and it is shown
that the corresponding matings define the same rational map $f$.}
\end{figure}

\textbf{Proof:} Denote the Airplane parameter by $q$ and fix
$n\in\{3,\,5,\,7,\,\dots\}$; let $c$ be the first center of period $n$
behind the Misiurewicz point $\gamma_\sM(5/12)$. The orbit of the
characteristic point $z_1$ is ordered as
\ban\nonumber
&&z_1<\gamma_c(5/12)<z_{n-1}<z_{n-3}<\dots<z_6<z_4<\alpha_c<\\[1mm]
&&<z_3<z_5<\dots<z_{n-2}<z_n<0<-\alpha_c<z_2 \ ; \label{shko} \ean
the critical orbit $(z_i^*)$ is similar with $z_n^*=0$. This
ordering is well-known from discussions of \v{S}harkovski\brevei{}
combinatorics. It can be checked with dynamic angles as follows: first, note
that the order of the critical orbit is compatible with the assumption that
$f_c\,:\,[z_1^*\,,\,0]\to[z_1^*\,,\,z_2^*]$ is strictly decreasing and
$f_c\,:\,[0,\,z_2^*]\to[z_1^*\,,\,z_3^*]$ is strictly increasing, so this
defines a unique real polynomial. Let $\Theta_1$ be the larger angle
at $z_1$ and denote its iterates under doubling by $\Theta_i$\,. Then
\ban \label{shka}
&&0<\Theta_2<1/6<\Theta_n<\Theta_{n-2}<\dots<\Theta_5<\Theta_3<1/3<\\[1mm]
&&<1/2<\Theta_1<7/12<\Theta_{n-1}<\Theta_{n-3}<\dots<\Theta_6<\Theta_4<2/3<1
\ , \nonumber \ean
since the derivative of the real polynomial $f_c(z)=z^2+c$ is negative for
$z<0$ and then $f_c$ swaps the lower and upper half-planes. Reading off binary
digits gives $\Theta_1=.\overline{10\,01\,01\dots01\,0}$, which is the largest
$n$-periodic angle less than $7/12=.10\,\overline{01}$. Reversing these
arguments, it follows that the center defined by $\gamma_\sM(\Theta_1)$ is real
and the orbit is as given by (\ref{shko}).
Each $z_i$ has two external angles, $\Theta_i$ and $1-\Theta_i$\,.
Note that $f_c:[0,\,\beta_c]\to[c,\,\beta_c]$ is increasing; taking preimages
of $z_i$ with respect to this branch gives strictly preperiodic points except
for $z_3$\,, which has a periodic preimage $z_2$ on the positive half-axis.

Now consider any parameter $p$ in the limb with rotation number $k/n$,
$k=(n-1)/2$. The wake is bounded by $0<\theta_-<\theta_+<1/3$. We have
$\theta_+=.\overline{01\,01\dots01\,0}$, since this is the largest
$n$-periodic angle less than $1/3=.\overline{01}$, or by sketching an
$n$-Rabbit. So $\theta_+=\Theta_3$ and $\theta_-=\Theta_5$\,; note that
$\Theta_1=1/2+\theta_+/4$ is an instance of the Douady Magic Formula from
Proposition~\ref{Pdmf}. --- The critical value $p$ of $f_p(z)=z^2+p$ is in
the sector at $\alpha_p$ bounded by the dynamic rays with the angles
$\theta_\pm$\,. This sector is mapped injectively for $n-1$ iterations; the
image sector contains $0$, $-\alpha_p$\,, and a unique point in
$f_p^{-1}(-\alpha_p)$\,. Thus the original sector contains unique preimages of
$\alpha_p$ with preperiods $n$ and $n+1$, respectively. Denote the angles of
the latter by $\theta_1<\dots<\theta_n$\,. Under $n$ iterations, these are
mapped to angles at $-\alpha_p$\,, such that $\theta_1$ gives the smallest
angle in $[1/2,\,1]$ and $\theta_n$ gives the largest angle in $[0,\,1/2]$. So
under $n+1$ iterations, $\theta_1$ is mapped to $\theta_+=\Theta_3$ and
$\theta_n$ is mapped to $\theta_-=\Theta_5$\,.

Next, let us look at the Airplane Julia set $\K_q$ with $Q(z)=f_q(z)=z^2+q$.
As the parameter was shifted from $c$ to $q$, the $n$-periodic points with
angles $\Theta_i$ moved holomorphically; in particular the pre-characteristic
points corresponding to $\pm z_n$ bound an interval containing the real slice
of the Airplane Fatou component around $0$. Consider the Fatou component of
$f_c$ at $z_3$\,; it defines an interval in $\K_q$\,, which contains a unique
preperiodic component $\Omega$ of preperiod $n-3$. Its largest antenna in the
upper halfplane has angles in a subset of
$[\Theta_5\,,\,\Theta_3]=[\theta_-\,,\,\theta_+]$. Since $f_q^{n-3}$ maps it
to the largest antenna on the upper side of the Fatou component around $0$,
$f_q^{n-2}$ maps it behind the component around $q$. Then it is behind the
component around $f_q(q)$, then to the right of the component at $0$, and
finally we see that $f_q^{n+1}$ maps the antenna of $\Omega$ to the interval
$(\gamma_q(4/7),\,\beta_q]$. Denote by $x_i$ the preimage of the $n$-periodic
point with angle $\Theta_i$\,, then $x_3$ has preperiod $n$ and the others
have preperiod $n+1$. On the other hand, the angles $\theta_i$ are the only
angles of preperiod $n+1$ in $(\theta_-\,,\,\theta_+)$ that are iterated to
some $\Theta_j$\,. Recalling that $\theta_1$ is iterated to
$\theta_+=\Theta_3$\,, we see that each $\theta_i$ with $i\neq1$ lands at some
$x_j$ with $j\neq3$. Denote the other angle by
$\theta_2'\,,\,\dots,\,\theta_n'$\.; it is in $(\theta_-\,,\,\theta_+)$ as
well, since the antenna is contained in an open half-strip bounded by these
rays and a real interval.

Finally, define the Misiurewicz parameters
$p_*=\gamma_\sM(\theta_1)=\dots=\gamma_\sM(\theta_n)$ and
$p_i=\gamma_\sM(\theta_i')$, $i=2,\,\dots,\,n$. Now $p_*$ is of $\alpha$-type
by construction, so it has preperiod $n+1$ and period $1$. The $p_i$ are
endpoints, since there is no other hyperbolic component of period $n$ in the
$k/n$-limb; they are pairwise different in particular. Note that for
$i=2,\,\dots,\,n$, the rays $\r_q(-\theta_i')$ and $\r_q(-\theta_i)$ land
together as well and the landing point never returns to this wake, so the two
rays are homotopic with respect to its orbit and to the real orbit of $q$, and
the precaptures are equivalent: by Proposition~\ref{Pshprec}, the matings
$Q\tmate P_i\eqg Q\tmate P_*$ agree, as do $P_i\tmate Q\eqg P_*\tmate Q$. ---
For the example of $k/n=2/5$, Figure~\ref{Fmult} shows the rays with angles
$-\theta_i\,,\,-\theta_i'$ landing pairwise at $\partial\K_q$, and the rays
with angles $\theta_i\,,\,\theta_i'$ landing at $\partial\K_{p_*}$\,, at a
preimage of $\alpha_{p_*}$ and at endpoints, respectively.

The landing pattern at $\partial\K_q$ is stable for parameters $q$ between
$c$ of period $n$ as above and the Airplane, but the relevant antenna will
bifurcate when $q$ is too far behind the Airplane.
\mybox

Note that we have constructed $n$ different matings giving the same rational
map, but in contrast to Theorem~\ref{Tchmun}, no upper bound on the
multiplicity is known in this case.
--- Assuming that the map $\M_{k/r}\to V_3$\,, $P\mapsto f\eqg P\tmate Q$ is
continuous, there will be self-intersections of the image corresponding to
these shared matings.

\subsection{Counterexamples to continuity of mating} \label{3c}
Geometric mating is not jointly continuous on the subset of $\M\times\M$ where
it can be defined. The first three examples below are due to Epstein
\cite{epstc, mateq}. Note that all of these techniques involve neutral
parameters, and that they do not exclude separate continuity. For specific
one-dimensional slices with $Q$ fixed, partial results on continuity have been
obtained by Dudko \cite{ddlam} and by Ma Liangang \cite{thesisma}.\\
\textit{--- Special thanks to Adam Epstein for explaining unpublished results.}
\begin{itemize}
\item Let $f_\lambda$ be a quadratic polynomial with a fixed point of
attracting multiplier $\lambda$. For $|\lambda|<1$, $|\mu|<1$ there are
explicit rational maps $F_{\lambda,\,\mu}\eqa f_\lambda\tmate f_\mu$\,.
Suppose $\lambda,\,\mu\to1$ tangentially, such that the third multiplier $\nu$
is constant. Then if $F_{\lambda,\,\mu}$ converges to a quadratic rational map,
it will depend on $\nu$, so there are oscillating sequences as well. Note that
convergence may depend on a normalization allowing the collision of the
respective fixed points; in a different normalization, $F_{\lambda,\,\mu}$
might converge to a map of degree one or to a constant as well.
\item Results on shared matings with cluster cycles by Sharland
\cite{shar1, shar2} are reported in Section~\ref{3s}. For rotation number
$1/n$, we have $f_n\eqg R_n\tmate Q_n\eqa P_n\tmate R_n$\,, where the
center parameters correspond to the following roots:
$r_n\sim\gamma_\sM(1/(2^n-1))=\gamma_\sM(2/(2^n-1))$,
$q_n\sim\gamma_\sM(-3/(2^n-1))=\gamma_\sM(-4/(2^n-1))$, and
$p_n\sim\gamma_\sM((2^{n-1}-1)/(2^n-1))=\gamma_\sM(2^{n-1}/(2^n-1))$.
Then $r_n\to r_0=1/4=\gamma_\sM(0)$, $q_n\to q_0=1/4=\gamma_\sM(0)$, and
$p_n\to p_0=-2=\gamma_\sM(1/2)$. Now if mating was continuous, we should have
$R_0\tmate Q_0\eqa P_0\tmate R_0$\,; both geometric matings exist, the former
has two parabolic basins and the latter has one. 
\item For a parabolic or bounded-type Siegel parameter $p$ on the boundary
of the main cardioid with angle $\theta$ and the real parameter $q$ defined
by the Douady Magic Formula $\Theta=1/2+\theta/4$ according to
Proposition~\ref{Pdmf}, consider the geometric mating $f_\theta\eqg P\tmate Q$,
which exists according tp Bl\'e-Valdez \cite{ble1, ble2}. When $\theta$ is
irrational, then $f_\theta^2(\infty)=0$, since the corresponding point in
$\K_q$ has the angles $\pm2\Theta=\pm\theta/2$ and the critical point of $P$
has $\theta/2$ as well. But when $\theta$ is rational, then either $0$ is
in a parabolic basin and $\infty$ is preperiodic, or there are disjoint cycles
of parabolic basins; in both cases $f_\theta^2(\infty)\neq0$. So approximating
a rational angle with irrational ones gives a contradiction to continuity.
\item Theorem~\ref{Tdiscbi} below uses similar ideas to show that the limit
is different from the expected one; since only rational angles are used, no
special arguments are needed to show matability. Here both $p_n$ and $q_n$
are Misiurewicz polynomials; a concrete example is given below as well.
\item Shared matings according to Theorem~\ref{Tchmun} can be used to
produce several counterexamples to continuity; here $p_n$ is hyperbolic and
$q_n$ is Misiurewicz. Again, the contradiction comes either from a different
number of parabolic Fatou cycles, or from an expected limit outside of the
Chebychev family. 
\item Different kinds of discontinuity may be expected in higher degrees.
E.g., with cubic polynomials $f_a(z)=z^3+az^2$, the mating
$f_a\tmate f_{-\overline a}$ gives an antipode-preserving rational map
\cite{bbm}. The former bifurcation locus shall be locally connected at
parabolic parameters, while the latter is not. So for suitable sequences of
postcritically finite polynomials, there will be an oscillatory behavior.
\end{itemize}

\begin{thm}[Discontinuity with bitransitive family]\label{Tdiscbi}
Consider a sequence of rational angles $\theta_n\to\theta_0$\,, such that
$\theta_n$ and $2\theta_n$ are preperiodic for $n\ge1$, $2\theta_0$
is periodic, and $\theta_0$ may be either unless $\theta_0$ and $2\theta_0$
belong to the same root. Set $p_n=\gamma_\sM(\theta_n)$ and
$q_n=\gamma_\sM(-2\theta_n)$ for $n\ge0$. Then the sequence of geometric
matings $f_n\eqg P_n\tmate Q_n$ does not converge to $f_0\eqg P_0\tmate Q_0$\,.
\end{thm}

\textbf{Proof:} First, note that $\theta$ and $2\theta$ are never in the same
limb, unless both are angles of the root. Thus all geometric matings under
consideration exist. Since the angle $\theta_n$ of $p_n\in\K_{p_n}$ is
complex conjugate to an angle $-\theta_n$ of $0\in\K_{q_n}$\,, there is a
direct ray connection between these two points, and the rational map
satisfies $f_n(0)=\infty$. We have $f_n\not\to f_0$ since $f_0(0)\neq\infty$:
while $z=\infty$ has an infinite orbit converging to a parabolic cycle of
$f_0$\,, $z=0$ either has a finite orbit or it converges to a different
parabolic cycle. --- This phenomenon seems to be analogous to parabolic
implosion, if we are looking at the polynomials $Q_n$ or at precaptures
according to Proposition~\ref{Phrc}: $q_n=\gamma_{q_n}(-2\theta_n)$ converges
to the critical value $q_0$ inside a parabolic Fatou component of $Q_0$\,, but
$\gamma_{q_0}(-2\theta_0)$ is a boundary point of this component. Of course,
parabolic implosion looks different for the rational maps here, since the
Julia set of $f_n$ is all of $\hat\C$.
\mybox

A concrete example is given by $\theta_n=u_n/2^{2n}$ with $u_n=(2^{2n-1}+1)/3$.
Then $p_n$ and $q_n$ are $\beta$-type Misiurewicz points, converging to the
Misiurewicz point $p_0=\i=\gamma_\sM(1/6)$ and the root
$q_0=-3/4=\gamma_\sM(1/3)$, respectively, and the matings do not converge to
the mating of the limits. Probably we have a parabolic $2$-cycle in both cases,
and Fatou components corresponding to a fat Basilica, but the limit of the
matings has $0$ and $\infty$ in different components of the Fatou set, while
the mating of the limits has $0$ in the Julia set at a preimage of the
parabolic fixed point.

\section{Short and long ray connections} \label{4}
We shall obtain explicit bounds on ray connections in Section~\ref{41},
discuss special irrational ray connections in Section~\ref{42}, search long
ray connections algorithmically in Section~\ref{43}, and give examples of
cyclic ray connections in Section~\ref{5}. The results provide partial
answers to Questions 3.1--3.3, 3.5--3.7, and 3.9 in \cite{mateq}.

\subsection{Bounding rational and irrational ray connections} \label{41}
When $p$ is postcritically finite, every biaccessible point $z\in\partial\K_p$
will be iterated to an arc $[-\beta_p\,,\,\beta_P]$, then to
$[\alpha_p\,,\,-\alpha_p]$, then to $[\alpha_p\,,\,p]$, and it stays within the
Hubbard tree $T_p\subset\K_p$\,. In \cite{m14}, Milnor discusses several
aspects of the geometric and the topological mating $P\tmate Q$ with
$p=q=\gamma_\sM(1/4)$. Every non-trivial ray connection
will be iterated to a connection between points on the Hubbard trees, since
every biaccessible point is iterated to the Hubbard tree $T_p$ or
$T_{\overline q}$\,. The two sides of the arcs of $T_p$ are mapped in a
certain way, described by a Markov graph with six vertices, such that only
specific sequences of binary digits are possible for external angles of
$T_p$\,. It turns out the only common angles of $T_p$ and $T_{\overline q}$ are
the $4$-cycle of $3/15$ and some of its preimages.
This fact implies that all ray connections between the Julia sets $\K_p$
and $\K_{\overline q}$ are arcs or trees of diameter at most $3$, so the
topological mating exists by Proposition~\ref{Palrcm}.

We shall consider an alternative argument, which is due to \cite{st3} in a
cubic situation. It gives weaker results in the example of $1/4\fmate1/4$, but
it is probably easier to apply to other cases: $T_q$ is obtained by cutting
away the open sector between the rays with angles $9/14$ and $11/14$, and its
countable family of preimages, from $\K_q$\,. So no $z\in T_{\overline q}$ has
an external angle in the open interval $(3/14,\,5/14)$, or in its preimages
$(3/28,\,5/28)$ and $(17/28,\,19/28)$. Now for every $z$ on the arc
$[\alpha_p\,,\,-\alpha_p]$, the angles on one side are forbidden. That shall
mean that the corresponding rays do not connect $z$ to a point in
$T_{\overline q}$\,, but to an endpoint of $\K_{\overline q}$ or to a
biaccessible point in a preimage of $T_{\overline q}$\,. This fact implies
that every ray-equivalence class has diameter at most four, which is weaker
than Milnor's result, but sufficient for the topological mating.

This argument shall be applied to another example, the mating of the
Kokopelli $P$ and the Airplane $Q$. Here $T_{\overline q}=T_q$ has no
external angle in $(6/7,\,1/7)$, and one side of $[\alpha_p\,,\,-\alpha_p]$
has external angles in $[1/14,\,1/7]$. Treating preimages of $\alpha_p$
separately, it follows that no other point in $\K_p$ is connected to two points
in $T_{\overline q}$\,, and we shall see that all ray-equivalence classes are
uniformly bounded trees. So the existence of the topological mating is
obtained without employing the techniques of Theorem~\ref{Trstrs}
by Thurston, Rees--Shishikura--Tan, and Rees--Shishikura. Moreover,
this approach works for geometrically finite and infinite polynomials as well.
E.g., $q$ may be any real parameter before the Airplane root, and $p$ be any
parameter in the small Kokopelli Mandelbrot set. 
Note however, that only the topological mating is obtained here, not the
geometric mating: there need not be a corresponding rational map.

To formulate the argument when $\K_q$ is locally connected but $Q$ is not
postcritically finite, we shall employ a generalized Hubbard tree $T_q$\,: it
is a compact, connected, full subset of $\K_q$\,, which is invariant under $Q$
and contains an arc $[\alpha_q\,,\,q]$. If $\K_q$ has empty interior and $q$ is
not an endpoint with irrational angle, there will be a minimal tree with these
properties. When $\K_q$ has non-empty interior, a forward-invariant topological
tree need not exist, but we may add closed Fatou components to suitable arcs to
define $T_q$\,. And when $q$ is an irrational endpoint, we shall assume that it
is renormalizable, and add complete small Julia sets to $T_{\overline q}$\,.
--- Note that in any case, every biaccessible point in $\K_q$ will be absorbed
by $T_q$\,, since $[\alpha_q\,,\,q]\subset T_q$\,.

\begin{prop}[Explicit bound on ray connections]\label{Pbdrcxpl}
Consider ray-equivalence classes for the formal mating $g=P\fmate Q$ of
$P(z)=z^2+p$ and $Q(z)=z^2+q$, with $\K_p$ and $\K_q$ locally connected, and
with a generalized Hubbard tree $T_q\subset\K_q$ as defined above. Now suppose
that there is an open set of angles, such that no external angle of
$T_{\overline q}$ is iterated to this forbidden set, and such that for an arc
$[\alpha_p\,,\,-\alpha_p]\subset\K_p$\,, the external angles on one side are
forbidden. Then:

$1$. Any point in $\K_p$ has at most one ray connecting it to a point in
the generalized Hubbard tree $T_{\overline q}$ of $\overline Q$\,.

$2$. All ray-equivalence classes have diameter bounded by eight, since each
class is iterated to a tree of diameter at most four.

$3$. Moreover, there are no cyclic ray connections, so the topological mating
$P\tmate Q$ exists according to Proposition~$\ref{Palrcm}$.
\end{prop}

\textbf{Proof:} 1. By assumption, $\alpha_p$ has at least one forbidden angle,
but there may be several allowed angles. Since these are permuted transitively
by iteration, none of them is connected to $T_{\overline q}$\,. In particular,
there is no ray connecting $\alpha_p$ to $\alpha_{\overline q}$\,, so $p$ and
$q$ are not in conjugate limbs. Suppose $z\in\partial\K_p$ is
not a preimage of $\alpha_p$\,. If it had two rays connecting it to points in
$T_{\overline q}$\,, this connection could be iterated homeomorphically until
both rays are on different sides of the arc $(\alpha_p\,,\,-\alpha_p)$,
contradicting the hypothesis since $T_{\overline q}$ is forward-invariant.
(Even if $z$ is precritical and reaches $0$ with both rays on one side, the
next iteration will be injective.)

2. Suppose $C$ is any bounded connected subset of a ray-equivalence class.
Iterate it forward (maybe not homeomorphically) until all of its preperiodic
points have become periodic, all critical and precritical points have become
postcritical, and all biaccessible points of $\K_{\overline q}$ have been
mapped into $T_{\overline q}$\,. So $C$ is a preimage of an eventual
configuration $C_\infty$, which is a subset of a ray-equivalence class of
diameter at most four, since it contains at most one biaccessible point of
$T_{\overline q}$\,. E.g., it might be a periodic branch point of $\K_p$
connected to several endpoints of $\K_{\overline q}$\,, or a point of
$T_{\overline q}$ connected to two or more biaccessible points of $\K_p$\,,
which are connected to endpoints of $\K_{\overline q}$ on the other side. In
general, taking preimages will give two disjoint sets of the same diameter in
each step, unless there is a critical value involved.

Now $C_\infty$ contains at most one postcritical point of $\K_{\overline q}$\,.
If there are several postcritical points of $\K_p$\,, then $C_\infty$ is
periodic, and preperiodic preimages contain at most one postcritical point of
$P$. So when pulling back $C_\infty$\,, the diameter is increased at most
twice, and it becomes at most 16. Actually, when $C_\infty$ has diameter $4$,
neither postcritical point can be an endpoint of $C_\infty$\,, and some
sketch shows that the diameter will become at most 8.

3. If $C$ is a cyclic ray connection, it will be iterated to a subset of a tree
$C_\infty$ according to item~2. This means that in the same step, both
critical points are connected in a loop $C'$, and $C''=g(C)$ is a simple arc
connecting the critical values $p\in\K_p$ and
${\overline q}\in\K_{\overline q}$\,. This cannot be a single ray, since
$p$ and $q$ are not in conjugate limbs. Suppose that $C''$ is of
the form $p-{\overline q}'-p'-{\overline q}$ with
${\overline q}'\notin T_{\overline q}$\,. Now ${\overline q}'$ is biaccessible,
so it will be iterated to $T_{\overline q}$\,, and then it must coincide with
an iterate of ${\overline q}$ by item~1. So $C''$ is not iterated
homeomorphically, and $p'$ must be critical or precritical. But then $C''$
would be contained in a finite periodic ray-equivalence class, and the critical
value of $P$ would be periodic, contradicting $p\in\partial\K_p$\,. The same
arguments work to exclude longer ray connections between the critical values
$p$ and ${\overline q}$.
\mybox

The following theorem provides large classes of examples. The parameter
$p$ is described by a kind of sector, and $q$ is located on some dyadic or
non-dyadic vein. More generally, $q$ may belong to a primitive or satellite
small Mandelbrot set, whose spine belongs to that vein. Let us say that $q$
is centered on the vein:

\begin{thm}[Examples of matings with bounded ray connections]\label{Tbdrcxpl}
When $p$ and $q$ are chosen as follows, with locally connected Julia sets,
the topological mating $P\tmate Q$ exists according to
Proposition~$\ref{Pbdrcxpl}$:

$a)$ The parameter $q$ is in the Airplane component or centered on the real
axis before the Airplane component, and $p$ in the limb $M_t$ with
rotation number $0<t\le1/3$ or $2/3\le t<1$.

$b)$ $q$ is centered on the non-dyadic vein to $\i=\gamma_\sM(1/6)$, and
$p\in\M_t$ with rotation number $0<t<1/2$ or $2/3<t<1$.

$c)$ $q$ is centered on the dyadic vein to $\gamma_\sM(1/4)$, and $p$ is
located between the non-dyadic veins to $\gamma_\sM(3/14)$ and
$\gamma_\sM(5/14)$. This means $p\in\M_t$ with $1/3<t<1/2$, or $p\in\M_{1/3}$
on the vein to $\gamma_\sM(3/14)$ or to the left of it, or $p\in\M_{1/2}$ on
the vein to $\gamma_\sM(5/14)$ or to the right of it. In particular, $p$ may
be on the vein to $\gamma_\sM(1/4)$, too.
\end{thm}

\textbf{Proof}: The case of $q$ in the main cardioid is neglected, because all
ray connections are trivial. We shall consider the angles of $\K_{\overline q}$
according to Figure~\ref{Fchsh}). When $Q$ has a topologically finite
Hubbard tree $T_q$\,, maximal forbidden intervals of
angles are found by noting that orbits entering $T_q$ must pass through
$-T_q$\,. See, e.g., Section~3.4 in \cite{core}. Denote the characteristic
angles of the limb $\M_t$ by $0<\theta_-<\theta_+<1$. For $p\in\M_t$\,, the arc
$[\alpha_p\,,\,\beta_p]$ has angles $\theta$ with $0\le\theta\le\theta_+/2$ on
the upper side and with $(\theta_-+1)/2\le\theta\le1$ on the lower side.

a) If $q$ is in the Airplane component or before it, the Hubbard tree is the
real interval $T_q=[q,\,q^2+q]$. If $q$ belongs to a small Mandelbrot set
centered before the Airplane, $T_q$ may contain all small Julia sets meeting
an arc from $q$ to $f_q(q)$ within $\K_q$\,.
Now no $z\in T_{\overline q}$ has an angle in $(6/7,\,1/7)$. So
Theorem~\ref{Tbdrcxpl} applies when $\theta_+/2<1/7$ or
$(\theta_-+1)/2>6/7$. The strict inequality is not satisfied for $t=1/3$ and
$t=2/3$. Then $\alpha_p$ and its preimages may be connected to three points in
the Hubbard tree of the Airplane, but the diameter is bounded by four as well.
Note that behind case a), with $q=\gamma_\sM(28/63)$ and $p=\gamma_\sM(13/63)$,
there is a ray connection of length six.

b) When $q$ is centered on the vein to $\gamma_\sM(1/6)$, the interval
$(11/14,\,1/14)$ is forbidden, so $(13/14,\,3/14)$ is forbidden for
$T_{\overline q}$\,. We need $\theta_+/2<3/14$ or
$(\theta_-+1)/2>13/14$.

c) For parameters $q$ centered on the vein to $\gamma_\sM(1/4)$, the interval
$(9/14,\,11/14)$ is forbidden, so $(3/14,\,5/14)$ is forbidden for
$T_{\overline q}$\,. We shall take its preimage
$(3/28,\,5/28)\cup(17/28,\,19/28)$ instead. When $p$ is between the veins to
$\gamma_\sM(3/14)$ and $\gamma_\sM(5/14)$, these two intervals are overlapping
in a sense: every $z\in(\alpha_p\,,\,-\alpha_p)$ has all angles on one side in
a forbidden interval. But then we have $p\in\M_t$ with $\theta_+/2<5/28$ or
$(\theta_-+1)/2>17/28$, so the forbidden intervals extend to $\pm\alpha_p$\,.
\mybox

\begin{xmp}[Bounded unlimited ray-equivalence classes]\label{Xreqmb}
Suppose $q$ is chosen according to item~a) or b), and $p$ is constructed as
follows. Take a primitive maximal component in the $1/3$-limb, then a primitive
maximal component in its $1/4$-sublimb, a primitive maximal component in its
$1/5$-sublimb \dots, then the limit $p$ has an infinite angled
internal address with unbounded denominators. $\K_p$ is locally connected by
the Yoccoz Theorem \cite{hy3, mlcj}, the topological mating exists according to
Theorem~\ref{Tbdrcxpl}, and there are branch points with any number of
branches. So ray-equivalence classes are bounded uniformly in diameter, but not
in size in the sense of cardinality.
\end{xmp}

\subsection{More on irrational ray connections} \label{42}
If two parameter rays with angles $\theta_-<\theta_+$ accumulate at the same
fiber of $\M$, it will intersect some dyadic vein in one point $c$, which is
called combinatorially biaccessible. $\K_c$ is locally connected and the
dynamic rays with angles $\theta_\pm$ land at the critical value $c$,
unless $c$ is parabolic. See the references in Section~4.4 of \cite{core}.
The following proposition shows that cyclic
ray connections for matings of biaccessible parameters can exist only in
special situations, since they cannot be preserved for postcritically finite
parameters behind them, where they are ruled out by Theorem~\ref{Trstrs}
of Rees--Shishikura--Tan. Compared to Proposition~\ref{Pbdrcxpl},
the situation is more general and the conclusion is weaker.

\begin{prop}[Cyclic irrational ray connections]\label{Pirrc1}
Consider the formal mating $g$ of $P(z)=z^2+p$ and $Q(z)=z^2+q$, with
parameters $p$ and $q$ not in conjugate limbs of $\M$.

a$)$ If $p$ is geometrically finite and $q$ is combinatorially biaccessible,
or vice versa, or both are geometrically finite, then $g$ does not have a
cyclic ray connection.

b$)$ If both $p$ and $q$ are combinatorially biaccessible and not
geometrically finite, then $g$ has a cyclic ray connection, if and only if
there is a ray connection between the critical values $p$ and $\overline q$.
\end{prop}

\textbf{Proof:} If both parameters are postcritically finite, the topological
mating exists according to Theorem~\ref{Trstrs}, and there
can be no cyclic ray connection by the Moore Theorem. For hyperbolic or
parabolic parameters, the ray connections will be the same as for the
corresponding centers. In general, a ray
connection between the critical values will have a cyclic preimage, so this
connection does not exist in case a). Conversely, a cyclic connection $C$ that
does not contain precritical points of the same generation, will give a
contradiction for postcritically finite parameters behind the current ones:
it may be iterated, possibly non-homeomorphically, to a cyclic connection
$C_\infty$ between points on the Hubbard trees, which are not critical or
precritical, and this connection $C_\infty$ would survive. To see this for
$P$, denote the external angles of the critical value $p$ by
$\theta_-<\theta_+$\,. Then no ray of $C_\infty$ will have an angle in
$(\theta_-/2\,,\,\theta_+/2)\cup((\theta_-+1)/2,\,(\theta_++1)/2)$. For
parameters $c$ behind $p$, the critical point is located in a strip
bounded by these four rays, so no precritical leaf can separate the rays
biaccessing points of $\K_p$ in $C_\infty$\,. (I have learned this technique
from Tan Lei.) The same argument applies to $\overline q$ and parameters behind
it. \mybox

The following proposition is motivated by Question~3.7 in \cite{mateq}. It
deals with angles $\theta$ that are rich in base $2$: the binary expansion
contains all finite blocks, or equivalently, the orbit of $\theta$ under
doubling is dense in $\R/\Z$. Angles with this property are rarely discussed
for quadratic dynamics, but they form a subset of full measure in fact.

\begin{prop}[Rich angles and irrational ray connections]\label{Pirrc2}
Suppose the angle $\theta$ is rich in base $2$.
Set $\theta_n=2^n\theta$ and $c_n=\gamma_\sM(\theta_n)$ for $n\ge1$. Then $c_n$
is a non-renormalizable endpoint of $\M$ with trivial fiber, $\K_{c_n}$ is a
dendrite, and the critical orbit is dense in $\K_{c_n}$\,.

$1$. For $n\neq m$ consider the formal mating $g$ of $P$ and $Q$, with $p=c_n$ and
$q=\overline{c_m}$\,. Then $g$ has a ray-equivalence class involving the
angle $\theta$, which is an arc of length four. $($Note that $n$ and $m$ may
be chosen such that $p$ and $q$ are not in conjugate limbs, but it is
unknown whether the topological or geometric mating exists.$)$

$2$. Let $\mathcal{X}_\theta\subset\M$ contain all parameters $c$, such that
$\theta$ is biaccessing $\K_c$\,. Then $\mathcal{X}_\theta$ is totally
disconnected, and it contains $c=-2$ and all $c_n$\,. So it has infinitely
many point components, and it is dense in $\partial\M$.
\end{prop}

\textbf{Proof:} Renormalizable and biaccessible parameters do not have dense
critical orbits. The orbit of an angle at the main cardioid is confined to a
half-circle \cite{bsoos}. By the Yoccoz Theorem \cite{hy3, mlcj},
$\K_{c_n}$ is locally connected with empty interior.

1. Assuming $n<m$, pull back the ray of angle $\theta_m$ connecting
postcritical points of $\K_p$ and $\K_{\overline q}$\,. This ray connects two
endpoints, so it forms a trivial ray-equivalence class. Since both points
are postcritical of different generations, the diameter is doubled twice
under iterated pullback (whenever there are two preimages, choose the component
containing an image of $\theta$).

2. For $c=-2$, every irrational angle is biaccessing, and for $c_n$\,, $\theta$
belongs to a critical or precritical point. By excluding all other cases,
$\mathcal{X}_\theta$ can contain only these and maybe other non-renormalizable,
postcritically infinite endpoints outside of the closed main cardioid, thus it
has only point components. So suppose that $\theta$ is biaccessing $\K_c$\,:\\
For a Siegel or Cremer polynomial of period 1, at most precritical points or
preimages of $\alpha_c$ are biaccessible \cite{szbc},
and the orbit of angles is not dense.\\
Pure satellite renormalizable parameters have only rational biaccessing angles
outside of the small Julia sets.\\
When $c$ is primitive renormalizable, the biaccessible points outside of the
small Julia sets are iterated to a set moving holomorphically with the
parameter, see Section~4.1 in \cite{core}. It is contained in a generalized
Hubbard tree $T_c$ in the sense of Proposition~\ref{Pbdrcxpl}.\\
When $c$ is postcritically finite or biaccessible, 
all biaccessible points are absorbed by a topologically finite tree $T_c$\,.
So their orbits are not dense in $K_c$ unless $T_c=\K_c$\,, which happens only
for $c=-2$.\\
It remains to show that $\mathcal{X}_\theta$ is dense in $\partial\M$: from a
normality argument it is known that $\beta$-type Misiurewicz points are dense.
For any Misiurewicz point $a=\gamma_\sM(\tilde\theta)$ there is a subsequence
with $\theta_n'\to\tilde\theta$\,. Then $c_n'\to a$, since
Misiurewicz points have trivial fibers \cite{sf2}. \mybox

\subsection{Searching long ray connections} \label{43}
Consider rational ray-equivalence classes for the formal mating $g=P\fmate Q$
with parameters $p,\,q$ in non-conjugate limbs of $\M$. A non-trivial periodic
ray connection requires pinching points in $\K_p$ and $\K_q$ with specific
angles, which exist if and only if the parameters $p,\,q$ are at or behind
certain primitive roots or satellite roots. So a longer ray connection means
that there are several relevant roots before the current parameters, and on the
same long vein in particular. Let us say that a ray connection is
\textbf{maximal}, if it is not part of a longer connection existing for
parameters behind the current ones. The following ideas were used to determine
all maximal ray connections algorithmically for ray periods up to $24$;
see Table~\ref{Ulist}.

\begin{table}[h!t!b!]
\begin{tabular}{|r||c|c|c|c|c|c|}\hline
 Per. & length $5$ & length $6$ & length $7$
 & length $8$ & length $10$ & length $12$\\ \hline\hline
$10$ & $32 \,+\, 0$ & $14 \,+\, 88$ & --- &
$0 \,+\, 2$ & --- & --- \\[1mm] \hline
$11$ & $76 \,+\, 0$ & $20 \,+\, 0$ & ---
 & --- & --- & --- \\[1mm] \hline
$12$ & $46 \,+\, 0$ & $24 \,+\, 264$ & ---
 & --- & --- & --- \\[1mm] \hline
$13$ & $226 \,+\, 0$ & $72 \,+\, 0$ & $2 \,+\, 0$
 & $2 \,+\, 0$ & --- & --- \\[1mm] \hline
$14$ & $285 \,+\, 0$ & $102 \,+\, 484$ & $4 \,+\, 0$
 & $0 \,+\, 14$ & $0 \,+\, 2$ & --- \\[1mm] \hline
$15$ & $540 \,+\, 0$ & $192 \,+\, 184$ & ---
 & --- & --- & --- \\[1mm] \hline
$16$ & $958 \,+\, 0$ & $338 \,+\, 1060$ & $4 \,+\, 0$
 & $2 \,+\, 10$ & $0 \,+\, 4$ & --- \\[1mm] \hline
$17$ & $1872 \,+\, 0$ & $584 \,+\, 0$ & $14 \,+\, 0$
 & $2 \,+\,  0$ & --- & --- \\[1mm] \hline
$18$ & $2814 \,+\, 0$ & $884 \,+\,2672$ & $22 \,+\, 0$
 & $6 \,+\, 24$ & $0 \,+\, 8$ & --- \\[1mm] \hline
$19$ & $5856 \,+\, 0$ & $1650 \,+\, 0$ & $26 \,+\, 0$
 & $6 \,+\, 0$ & --- & --- \\[1mm] \hline
$20$ & $9534 \,+\, 0$ & $2890 \,+\,5244$ & $58 \,+\, 0$
 & $4 \,+\, 42$ & $0 \,+\, 8$ & --- \\[1mm] \hline
$21$ & $16978 \,+\, 0$ & $4900 \,+\, 898$ & $64 \,+\, 0$
 & $4 \,+\, 0$ & --- & --- \\[1mm] \hline
$22$ & $30180 \,+\, 0$ & $8423 \,+\, 10928$ & $126 \,+\, 0$
 & $18 \,+\, 132$ & $0 \,+\, 20$ & $0 \,+\, 2$ \\[1mm] \hline
$23$ & $55676 \,+\, 0$ & $15300 \,+\, 0$ & $172 \,+\, 0$
 & $18 \,+\, 0$ & --- & --- \\[1mm] \hline
$24$ & $95830 \,+\, 0$ & $25968 \,+\, 25312$ & $242 \,+\, 0$
 & $24 \,+\, 96$ & $0 \,+\, 28$ & --- \\[1mm] \hline
\end{tabular}
\caption[]{\label{Ulist}%
The length of maximal periodic ray connections depending on the ray period.
The first number counts unordered pairs of periodic parameters with
primitive-only connections, the second number is the connections including a
satellite cycle. Length $\le4$ is ubiquitous, 
length $5$ appears already for periods $7$ and $9$, while length $6$ happens
for periods $4$ and $6$--$9$ as well. Length $9$ and $11$ was not found for
periods $\le24$.}
\end{table} 
\begin{itemize}
\item Suppose $\r(\theta_1)$--$z_p$--$\r(\theta_2)$--$z_q$--$\r(\theta_3)$ is
a step in the ray connection, then $\theta_1$ and $\theta_2$ belong to a
cycle of angle pairs for $\K_p$\,, so there is a root before $p$ with
characteristic angles iterated to $\theta_1$ and $\theta_2$\,. Likewise,
there is a root before $\overline q$, whose characteristic angles are
iterated to $\theta_2$ and $\theta_3$\,. Conversely, given the angles
$\theta_\pm$ of a root before $p$, we may determine conjugate angles for
iterates of $\theta_+$ under doubling, and check whether the root given by an
angle pair is before $\overline q$; it is discarded otherwise. So we record
only the angle pairs of roots, and forget about the number of iterations and
about which class in a cycle contains which characteristic point. Note that
there is an effective algorithm to determine conjugate angles
\cite{bks, renorm}, probably due to Thurston.
\item A maximal ray connection should be labeled by highest relevant roots on
the respective veins. However, a brute-force search starting with these roots
will be impractical: varying both $p$ and $\overline q$ independently is too
slow, and searching $\overline q$ depending on $p$ requires to match different
combinatorics on two sides, since the characteristic point $z_p$ corresponding
to the highest root may be anywhere in the ray-equivalence class. So the idea
is to run over all roots $p_1$\,, try to build a maximal ray connection on one
side of the corresponding characteristic point, and to quit if the connection
can be continued on the other side of that point.
\item When a pinching point of satellite type is reached under the recursive
application of the conjugate angle algorithm, we may double the length and
stop. Alternatively, two separate algorithms may be used, one finding
primitive-only ray connections starting from the first pinching point, and
another one starting with the satellite-type point in the middle of the
periodic ray-equivalence class.
\end{itemize}
For period $22$, this algorithm has recovered the example given to Adam
Epstein by Stuart Price \cite{mateq}: for $p$ behind
$\{1955623/4194303,\,1955624/4194303\}$ and $q$ behind
$\{882259/4194303,\,882276/4194303\}$ there is a periodic ray-equivalence class
of diameter $12$. For $1/2$-satellites only, the same algorithm was used for
periods up to $40$ in addition; this produced another example of diameter $14$
for period $32$, with $p$ behind
$\{918089177/4294967295,\,918089186/4294967295\}$ and $q$ behind
$\{1998920775/4294967295,\,1998920776/4294967295\}$. Note that, e.g., taking
$p$ and $q$ as the corresponding centers, the formal mating will have
non-postcritical long ray connections and the geometric mating shows
clustering of Fatou components. For suitable preperiodic parameters behind
these roots, the formal mating has long periodic ray-equivalence classes with
postcritical points from both orbits, and preperiodic classes may have
twice or up to four times the diameter of the periodic classes.
--- There are several open questions on long ray connections:
\begin{itemize}
\item What are possible relations between the linear order of roots on the
veins to $p$ and $q$, and the order of pinching points within a ray-equivalence
class?
\item For the lowest period with a particular diameter of a ray-equivalence
class, is there always a $1/2$-satellite involved? 
\item Is there a whole sequence with similar combinatorics and increasing
diameters? If it converges, does the limit show non-uniformly bounded ray
connections? 
Does the geometric mating of the limits exist?
If not, does it have infinite irrational ray connections?
\item Are there only short ray connections for self-matings and for matings
between dyadic veins of small denominator, even though the Hausdorff dimension
of biaccessing angles is relatively high according to \cite{dsj}?
\end{itemize}

\section{Cyclic ray connections} \label{5}
\begin{figure}[h!t!b!]
\unitlength 0.001\textwidth 
\begin{picture}(990, 420)
\put(10, 0){\includegraphics[width=0.42\textwidth]{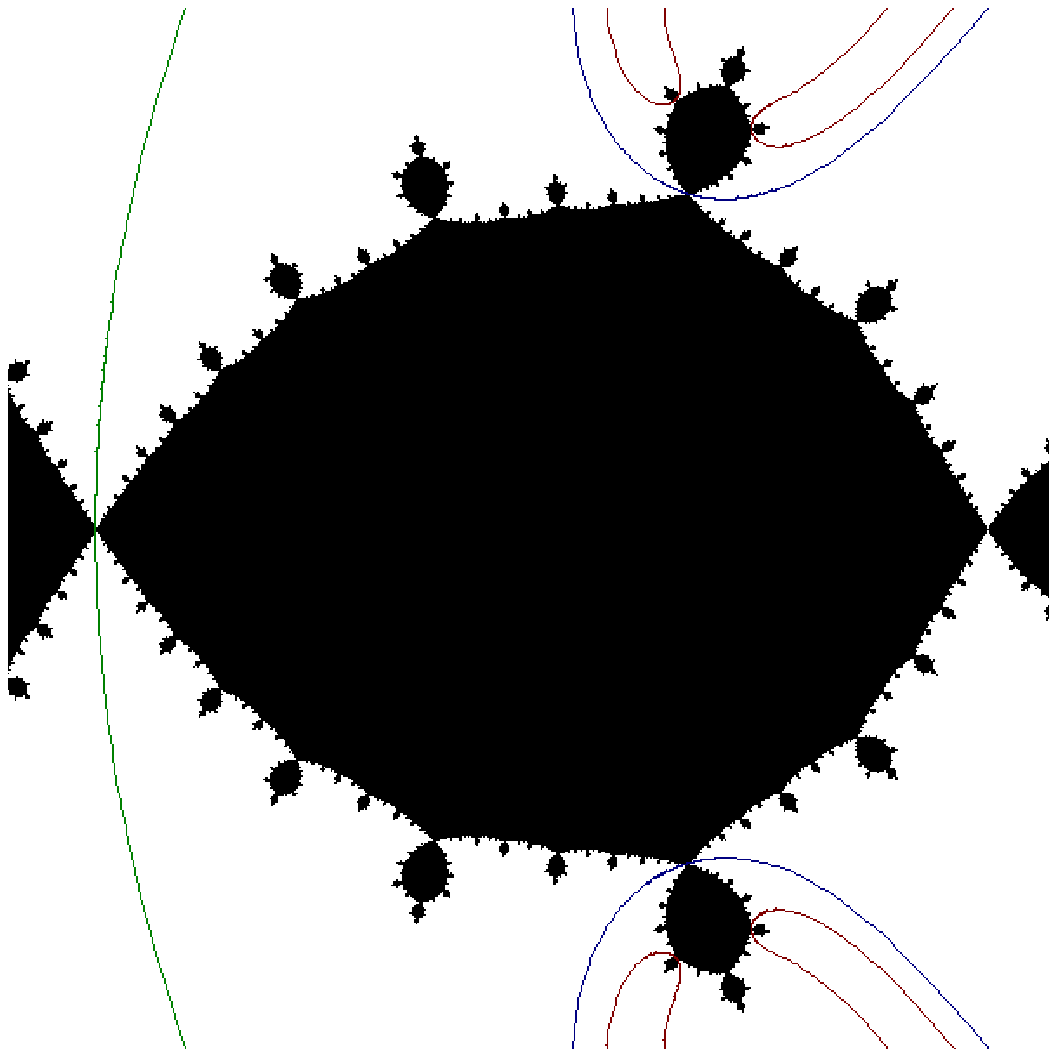}}
\put(570, 0){\includegraphics[width=0.42\textwidth]{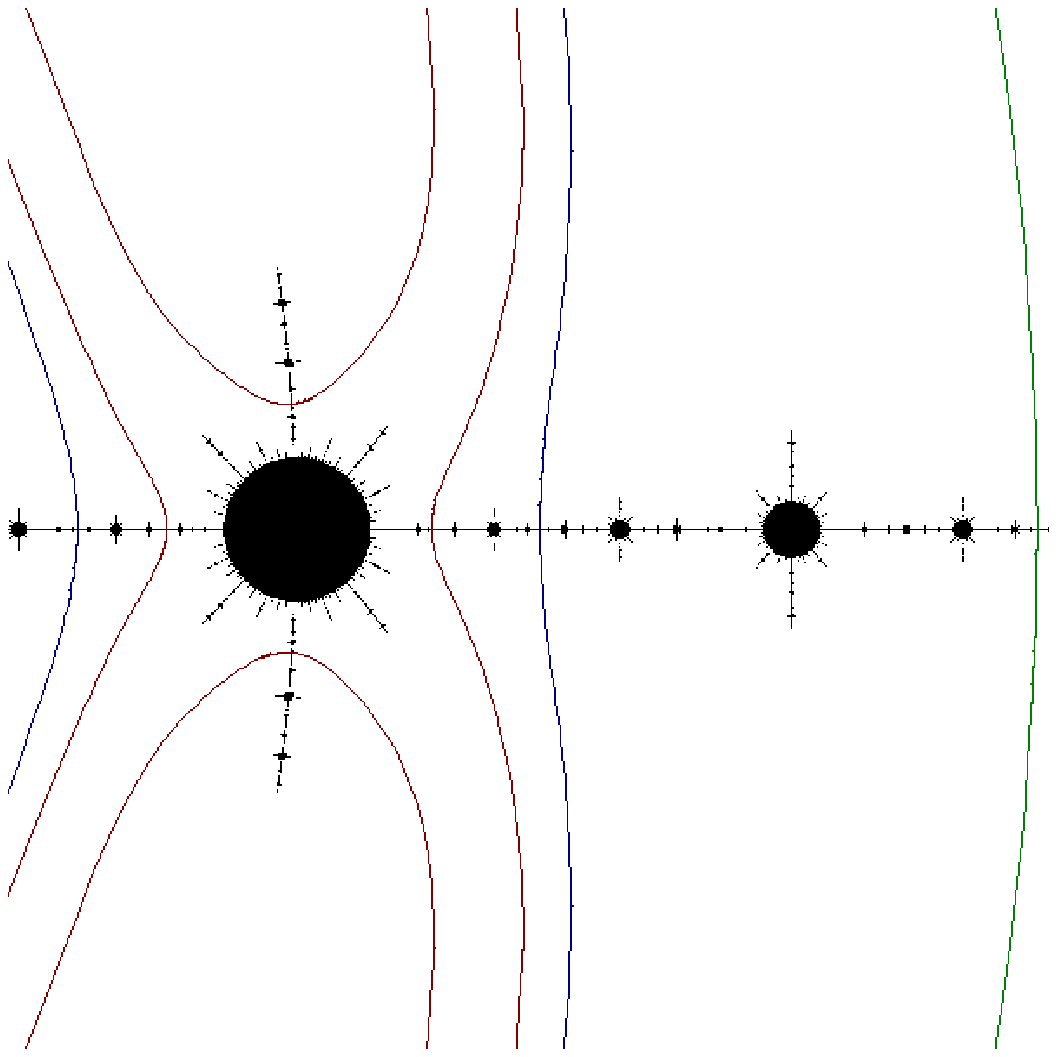}}
\thinlines
\multiput(10, 0)(560, 0){2}{\line(1, 0){420}}
\multiput(10, 420)(560, 0){2}{\line(1, 0){420}}
\multiput(10, 0)(560, 0){2}{\line(0, 1){420}}
\multiput(430, 0)(560, 0){2}{\line(0, 1){420}}
\end{picture} \caption[]{\label{FrenAB}
The formal mating $g$ of the Airplane $\K_p$ (on the right) and the Basilica
$\K_q$ (shown rotated on the left). The green ray connection $C_0$ has the
angles $5/12$ and $7/12$. Suitable preimages $C_1$ (blue), $C_2$ (red), \dots
form nested cycles around the critical value component of the Airplane. The
nested domains are typical of primitive renormalization. The canonical
obstruction of $g$ is discussed in Figure~2 of \cite{qmate}.}
\end{figure}

First we shall construct cyclic ray connections for the formal mating $g$ of
the Airplane $P(z)=z^2-1.754877666$ and the Basilica $Q(z)=z^2-1$. See
Figure~\ref{FrenAB}. All biaccessing rays of $Q$ are iterated to the
angles $1/3$ and $2/3$ at $\alpha_q=\alpha_{\overline q}$\,. Denote by $C_0$
the cyclic connection formed by the rays with angles $5/12$ and $7/12$. Pulling
it back along the critical orbit of the Airplane gives nested cycles $C_n$
around the critical value $p$, since $g^3$ is proper of degree $2$
from the interior of $C_1$ to the interior of $C_0$\,. Now $C_n$ has $2^n$
points of intersection with $\K_p$\,, so its length is not uniformly bounded
as $n\to\infty$. Moreover, $C_n$ connects points $x_n$ converging to
$x_\infty=\gamma_p(3/7)=\gamma_p(4/7)$ to points $x_n'$ converging
to $x_\infty'=\gamma_p(25/56)=\gamma_p(31/56)$. But these four rays are landing
at endpoints of the Basilica, so the landing points $x_\infty\neq x_\infty'$
on the Airplane critical value component are not in the same ray-equivalence
class. Thus the ray-equivalence relation is not closed. In fact, the limit set
of $C_n$ contains the boundary of the Fatou component, which meets uncountably
many ray-equivalence classes. I am not sure what the smallest closed
equivalence relation, or the corresponding largest Hausdorff space, will look
like: it shall be some non-spherical quotient of the Basilica, with a countable
family of simple spheres attached at unique points. This Hausdorff obstruction
has been obtained independently by Bartholdi--Dudko [private communication].
--- More generally, we have:

\begin{thm}[Unbounded cyclic ray connections]\label{Tucrc}
Suppose $p$ is primitive renormalizable of period $m$ and $\K_p$ is locally
connected. Then there are parameters $c_*\prec c_0\prec p$\,, such that for
all parameters $q$ with $\overline{q}$ on the open arc from $c_*$ to $c_0$\,,
the formal mating $g=P\fmate Q$ has non-uniformly bounded cyclic ray
connections. Moreover, these are nested such that the ray-equivalence
relation is not closed. So the topological mating $P\tmate Q$ is not defined
on a Hausdorff space.
\end{thm}

\textbf{Proof:} In the dynamic plane of $\K_p$\,, denote the small Julia set
around the critical value $p$ by $\K_p^m$\,. There are preperiodic
pinching points with
$\alpha_p\preceq x_\ast\prec x_0\prec x_1\prec\K_p^m\prec x_1'$\,, such that
$P^m$ is a $2$-to-$1$ map from the strip between $x_1$ and $x_1'$ to the wake
of $x_0$\,. Restricting these sets by equipotential lines in addition, we
obtain a polynomial-like map, which is a renormalization of $P$\,. If
the pinching points are branch points, the bounding rays must be chosen
appropriately. We assume that $x_1$ and $x_1'$ are iterated to $x_0$ but never
behind it, and $x_0$ is iterated to $x_*$ but never behind it. More generally,
$x_*$ may be a periodic point. The construction of these points is
well-known from primitive renormalization; see \cite{sf2, tlp, renorm}.

Since the points $x_*$ and $x_0$ are characteristic in $\K_p$\,,
there are corresponding Misiurewicz points $c_*$ and $c_0$ in $\M$.
(If $x_*$ is periodic, then $c_*$ is a root.) When the parameter $\overline{q}$
is in the wake of $c_*$\,, or in the appropriate subwake, then $x_0$ will be
moving holomorphically with the parameter and keep its external angles.
When $\overline{q}$ is chosen on the regulated arc from $c_*$ to $c_0$\,, then
$\K_{\overline q}$ will be locally connected. In $\K_{\overline q}$ the point
corresponding to $x_0$ has the same external angles as in $\K_p$\,, and no
postcritical point is at this point or behind it.
Thus the four rays defining the strip between $x_1\,,\,x_1'\in\K_p$ are
landing in a different pattern at $\K_{\overline q}$\,.

Now consider the formal mating $g=P\fmate Q$. We shall keep the notation
$x_i\,,\,p,\,\K_p\,,\,\overline{q},\,\K_{\overline q}$ for the corresponding
points and sets on the sphere. Since the two rays bounding the wake of $x_0$\,,
or the relevant subwake, are landing together at $\K_{\overline q}$\,, they
form a closed ray connection $C_0$\,. Its preimage is a single curve consisting
of four rays, two pinching points in $\K_p$\,, and two pinching points in
$\K_ {\overline q}$\,. This can be seen on the sphere, since $C_0$ is
separating the critical values of $g$, or in the dynamic plane of
$\overline{q}$, since $\overline{q}$ is not behind the point corresponding to
$x_0$\,. Now the new curve is pulled back with $g^{m-1}$ to obtain $C_1$\,,
which is a closed curve connecting $x_1$ and $x_1'$ to two pinching points in
$\K_{\overline q}$\,. By construction, $g^m$ is proper of degree $2$ from the
interior of $C_1$ to the interior of $C_0$\,, and the former is compactly
contained in the latter. $g^m$ behaves as a quadratic-like map around
$\K_p^m$\,, but only points below the equator will converge to the small Julia
set under iterated pullback.

Define the curves $C_n$ inductively; they form strictly nested closed curves 
and the number of rays is doubled in each step. E.g.,
$C_2$ is intersecting $\K_p$ in four points. The two preimages $x_2$ and $x_2'$
of $x_1$ are located between $x_1$ and $x_1'$\,, while the two preimages of
$x_1'$ belong to decorations of $\K_p^m$ attached at the points with
renormalized angles $1/4$ and $3/4$. We have
$x_0\prec x_1\prec x_2\prec\dots\prec\K_p\prec\dots\prec x_2'\prec x_1'$\,.
The limits $x_\infty$ and $x_\infty'$ are the small $\beta$-fixed point of
$\K_p^m$ and its preimage, the small $-\beta$. Now $x_n$ and $x_n'$ are
connected by $C_n$\,, but $x_\infty$ and $x_\infty'$ are not ray-equivalent,
because the former is periodic and the latter is preperiodic. \mybox

More generally, $\overline q$ may be any parameter in the strip between
$c_*$ and $c_0$\,, as long as its critical orbit does not meet the point
corresponding to $x_0$ or get behind it. 
--- Note that by taking iterated preimages of a finite ray-equivalence tree,
you will merely get uniformly bounded trees: the diameter can be increased only
when a critical value is pulled back to a critical point, which can happen at
most twice according to Proposition~\ref{Psrec}: a finite irrational tree
cannot be periodic, so it does not contain more than one postcritical point
from each polynomial.

The program Mandel provides several interactive features related to the
Thurston Algorithm. It is available from
\href{http://www.mndynamics.com/indexp.html}{www.mndynamics.com}\,.
A console-based implementation of slow mating is distributed with the
\href{http://arxiv.org/abs/1706.04177}{preprint} of \cite{qmate}.
\end{document}